\documentclass[11pt]{article}   
\usepackage{amssymb,amscd,latexsym} 
\usepackage{graphicx,color}  
\usepackage{amsmath}
\usepackage{amsthm}
\usepackage{mathdots}
\usepackage[all]{xy}
\newcommand{\rar}{\rightarrow}
\newcommand{\lar}{\longrightarrow}
\newcommand{\llar}{-\kern-5pt-\kern-5pt\longrightarrow}
\newcommand{\surjects}{\twoheadrightarrow}

\newtheorem{Theorem}{Theorem}[section]
\newtheorem{Lemma}[Theorem]{Lemma}
\newtheorem{Corollary}[Theorem]{Corollary}
\newtheorem{Proposition}[Theorem]{Proposition}
\newtheorem{Remark}[Theorem]{Remark}

\newtheorem{Conjecture}[Theorem]{Conjecture}
\newtheorem{Definition}[Theorem]{Definition}
\newtheorem{Question}[Theorem]{Question}
\def\sqr#1#2{{\vcenter{\hrule height.#2pt
        \hbox{\vrule width.#2pt height#1pt \kern#1pt
            \vrule width.#2pt}
        \hrule height.#2pt}}}
\def\phi{\varphi}
\def\demo{\noindent{\bf Proof. }}
\def\square{\mathchoice\sqr64\sqr64\sqr{4}3\sqr{3}3}
\def\qed{\hspace*{\fill} $\square$}

\def\xx{{\bf x}}

\def\tt{{\bf t}}

\def\fm{{\mathfrak m}}

\def\hht{{\rm ht}\,}

\def\ker{{\rm ker}\,}

\def\restr{{\kern-1pt\restriction\kern-1pt}}

\def\pp{{\mathbb P}}

\begin{document}
\begin{center}
{\Large{\bf\sc Linearly presented perfect ideals of codimension $2$ in three variables}}
\footnotetext{AMS Mathematics
	Subject Classification (2010   Revision). Primary 13A30, 13C40, 13D02, 13H10; Secondary  13A02,14M12.}

\vspace{0.3in}

{\large\sc Andr\'e Doria}\footnote{Under a CAPES-PNPD post-doctoral fellowship} \quad\quad
{\large\sc Zaqueu Ramos} \quad\quad
 {\large\sc Aron  Simis}\footnote{Partially
 	supported by a CNPq grant (302298/2014-2), by a CAPES-PVNS Fellowship (5742201241/2016). This author thanks the ICMC of the University of S\~ao Paulo, S\~ao Carlos, Brazil, where part of this work has been done, under a Senior Visitor Fellowship.}

\end{center}


\begin{abstract}
The goal of this paper  is the fine structure of the ideals in the title, with emphasis on the properties of the associated Rees algebra and the special fiber. The watershed between the present approach and some of the previous work in the literature is that here one does not assume that the ideals in question satisfy the common generic properties. One exception is a recent work of N. P. H. Lan which inspired the present work. Here we recover and extend his work.
We strongly focus on the behavior of the ideals of minors of the corresponding so-called Hilbert--Burch matrix and on conjugation features of the latter.  We apply the results to three important models: linearly presented ideals of plane fat points, reciprocal ideals of hyperplane arrangements and linearly presented monomial ideals.
\end{abstract}

\section*{Introduction}

The existing literature on codimension $2$ perfect ideals is so vast that yet a new paper on the subject would look like a temerity, and yet if one forgoes the usual generic conditions, many questions are largely open.
The main purpose of this work is to pursue some of these questions, with a focus on the associated Rees algebra of the ideal.

Let $R$ denote a commutative Noetherian ring  and let $I\subset R$ stand for an ideal. 
Let $ \mathcal{R}_R(I):=R[It]\subset R[t]$ denote the Rees  algebra of $I$, where $t$ is an indeterminate over $R$. 
Recall that there is a structural graded $R$-algebra surjection $\mathcal{S}_R(I)\twoheadrightarrow \mathcal{R}_R(I)$, where $\mathcal{S}_R(I)$ denotes the symmetric algebra of $I$. For some of the early findings about this map and its related ideal theory and homology one may refer to \cite{Trento}, \cite{Wolmbook1} and the references there.
One such ideal theoretic notion  is the following:  the ideal $I$ is of  {\em linear type} if this map is injective.
A major property of such an ideal is the known restriction on its number of generators, namely, that the latter is locally everywhere bounded by the dimension of the ambient ring.

A whole lot of published work in commutative algebra has a focus on how to approach ideals which are not of linear type -- and, in particular, do not satisfy the above restriction on the number of generators  -- but still share a few common features that may be reflected in the behavior of the associated algebras.

Surprisingly, even for a homogeneous perfect ideal $I$ of codimension $2$ in a standard graded polynomial ring $R=k[\xx]=k[x_1,\ldots, x_d]$ over a field $k$, which is moreover generated in on single degree, there is no encompassing theory in all cases. In  addition to the Rees algebra, one is naturally led to study the so-called special fiber algebra. One relevant aspect of the latter is that it isomorphic in this case to the homogeneous defining  ring of the image of the rational map defined by the linear system of the generators of $I$.
There are satisfactory results in the case where $I$ satisfies the so-called $(G_d)$ condition, a property regarding the number of generators of $I$ locally off the irrelevant ideal $(\xx)$. This condition has a hidden impact on the codimension of an ideal of minors of lower size of the Hilbert--Burch matrix, making them sort of generic.
Thus, if $d=3$ the property means that $I$ is generically a complete intersection.
We refer to \cite{MU} and the references thereof for the main results along this line.

The case $d=2$ has been largely dealt with by many authors, mainly envisaging elimination theory (see \cite{Cox_et_al_1}, \cite{Cox_et_al_2},\cite{syl1}, \cite{Kustin_et_al}). In this environment one takes huge profit from the fact that the ideal $I$ is primary to the maximal irrelevant ideal. In any case, the methods and ideas, if not some of the intermediate results, tend to be quite apart from the ones in this work.

One first attempt to get away from the $(G_d)$ condition or others of similar nature is a recent result of N. Lan (\cite{Lan}), assuming that $d=3$ and that $I\subset k[x,y,z]$ is a linearly presented codimension $2$ perfect ideal.
In this work, as in \cite{Lan}, we only deal with the case where $d=3$ and assume throughout that $I$ is linearly presented. A reason for the latter is twofold: first, the technology becomes easier and second, one can expect that the Rees algebra itself be a Cohen--Macaulay ring in interesting cases.
In fact, allowing for minimal syzygies of higher degree usually disrupts the Cohen--Macaulayness of the Rees algebra (see \cite[Introduction]{Kustin_et_al} for $d=2$ and \cite[Section 4.1]{SymbCrem} for $d=3$).

With Lan we consider the $(n+1)\times n$ syzygy matrix $\phi$ of the ideal $I$ (i.e., its Hilbert--Burch matrix) and look at the nested sequence of its ideals of minors of various sizes:
$$I=I_n(\phi)\subset I_{n-1}(\phi)\subset\cdots\subset I_1(\phi),$$
where $I_1(\phi)$ is assumed to have codimension $3$ (a harmless hypothesis as the matrix has only linear entries).
Whereas Lan's main assumption is that $I_2(\phi)$ has codimension $2$, we impose no such restriction at the outset and instead look at the first index $t$ such that $I_{t+1}$ has codimension $2$. For lack of better terminology we call this index the {\em chaos invariant} of $\phi$ (or of $I$, keeping in mind the relevant fact that it is invariant under conjugation).
The chaos invariant will play a major role in this work.

We now briefly describe the contents of the sections and at the end state the main results of the paper. Let the reader have in mind that the main goal is the structure and properties of both the Rees algebra and the special fiber, such as Cohen--Macaulayness and normality, and the fiber type condition.

Section 1 is a brief reminder about the main ideal theoretic gadgets.
Most of the used terminology can be found here or in the suggested references thereof.

Section 2 develops the significant properties of the chaos invariant as regards the structure of the Hilbert--Burch matrix and the associated rational map defined by the linear system of the minors. It also captures the intervention of the reduction number in its relation to the Rees algebra and the special fiber. The case of reduction number $2$ is specially dealt with.

Section 3 calls upon a second main character, namely, a certain submatrix of the so-called Jacobian dual matrix developed in various sources (\cite{AHA}, \cite{MU}, \cite{SymBir}) and originally introduced in \cite{dual}. This submatrix undergoes several stages of virtual cloning throughout the text, all satisfying the property of being $1$-generic in the sense of \cite{Eisenbud2} and \cite[Chapter 9]{Harrisbook}. Geometrically, it will be responsible for taking us inside the environment of rational normal scrolls, hence yielding an upper hand on the properties of Cohen--Macaulayness and normality.
It is in this section that we are able to recover the results of \cite{Lan} by different methods, perhaps in a slightly simpler or more conceptual way.

Finally, Section 4 is devoted to the detailed study of a triplet of classes of these ideals: linearly presented ideals of plane fat points, reciprocal ideals of hyperplane arrangements and linearly presented monomial ideals.

The first class has of course been vastly studied by many experts. In order to apply to our scheme, one ought to understand as a preliminary step when the ideal is generated in one single degree -- linear presentation coming next as an even more difficult question.
Thus, here we limit ourselves to a very special class of fat ideals as considered in \cite{fat1}, \cite{fat2}). Some of these allow by a plane quadratic transformation to switch to a linearly presented codimension $2$ perfect ideal.
Unfortunately, as of now, this is the size of the possible illustration we can accomplish within this class of ideals.

The second class, subsumed in the first one, comes in the way of illustrating how some of the results of \cite{blup} happen in dimension $3$ in the light of the present approach. Of course, quite a bit about this class is well-known to the experts (see in particular\cite{Schenck}, \cite{Te}, \cite{PrSp}, \cite{blup}). Here we give an encore in dimension $3$ of the result of \cite{blup} that the Rees algebra is Cohen--Macaulay, as in this low dimension the argument becomes more pliable.

The class of codimension $2$ linearly presented perfect ideals generated by monomials is a surprisingly intricate world even if it appears very restrictive.
Alas, we could for the moment only treat adequately the case where the ideal has only two minimal primes (leaving out the case of three minimal primes).
An undoubtedly unexpected intervention has been the normality of binary monomial ideals via \cite{GSWR}.
Moreover, complete results about the Rees algebra and the special fiber were only accomplished in two special subclasses of these ideals.
Fortunately, here we have been fairly complete proving the Cohen--Macaulayness of the envisaged algebras, plus normality and/or fiber type property. 
This leaves huge ground unaccounted for at this moment. 

On a fingertip list, the main results of the paper are as follows:
\begin{itemize}
	\item 
 Theorem~\ref{Fitting_main} shows that the generators of $I$ define a birational map of $\pp^2$ onto its image in $\pp^{n-1}$, with inverse given by by a linear system of quadrics; in particular, $I$ has maximal analytic spread and depth$(R/I^2)=0$.

\item Theorem~\ref{red_no2} says that if the reduction number of $I$ is at most $1$ locally at the associated primes of $R/I$ then the actual reduction number of $I$ is at most $2$ if and only if the Rees algebra of $I$ is Cohen--Macaulay if and only if the special fiber of $I$ is Cohen--Macaulay.

\item Theorem~\ref{quadrics_is_scroll} proves that if $u$ denotes the chaos invariant of $I$, then a suitable $2\times (n-u)$ submatrix of the Jacobian linear dual of $I$ is $1$-generic; in particular its $2$-minors define an $(u+1)$-dimensional rational normal scroll in $\pp^{n-1}$, hence an arithmetically Cohen--Macaulay and projectively normal variety.

\item Theorem~\ref{Lan1} recovers the result of Lan's with a few additional contents, namely,  if the chaos invariant is $1$ (minimal possible) then (i) all the ideals $I_t(\phi)$, for $t\neq 1, n-1$, have one and the same minimal prime $q$, and locally at any associated prime other than $q$, $I$ is a complete intersection; (ii) the Rees algebra of $I$ is normal, Cohen--Macaulay and of fiber type.

For the last two theorems, one assumes that $I$ is minimally generated by monomials and that $R/I$ has exactly two minimal primes. For such an ideal, the Hilbert--Burch matrix has the following form up to conjugation:
\begin{equation}\label{canonical_monomial}
\varphi=\left(\begin{array}{ccccccccccc}
z&0&0&\ldots&0&0&0\\
-c_1&z&0&\ldots&0&0&0\\
0&-c_2&z&\ldots&0&0&0\\
\vdots&\vdots&\vdots&\ddots&\vdots&\vdots&\vdots\\
0&0&0&\ldots&z&0&0\\
0&0&0&\ldots&-c_{n-3}&z&0\\
0&0&0&\ldots&0&-c_{n-2}&z\\
0&0&0&\ldots&0&0&-c_{n-1}
\end{array}\right)
\end{equation}
where $c_i\in\{x,y\}$ for each $1\leq i\leq n-1.$
Call the sequence of entries $c_1 \, c_2\,\ldots\, c_{n-1}$ the basic entry sequence of $I$.

\item Theorem~\ref{alternating} tells that if the basic sequence has the shape $x\,y\,x\,y\,x\,y\,\ldots$, then the Rees algebra and the special fiber of $I$ are both normal and Cohen--Macaulay.
 
\item Theorem~\ref{separating} shows that if the basic sequence has the shape $x\,x\,x\,\ldots \,x\,y\,y\,\ldots\, y$, then the Rees algebra and the special fiber of $I$ are both Cohen--Macaulay, and the Rees algebra is of fiber type.
\end{itemize}

A word against excessive optimism: although our effort throughout has been to look for a meaningful situation in which the Rees algebra and the special fiber are Cohen--Macaulay (possibly, normal too) and the ideal is of fiber type, -- luckily enough, this includes all ideals with chaos invariant $1$ -- it is quite possible that the full majority of the ideals under the title fail to have such niceties.
The first author has developed many of lines of computation to have a feel for how often these niceties occur and the overall feeling is that of sparseness.

The third author heartily thanks N. Lan for telling about his work during a conversation in a lovely bistro in Nice.

\section{Ideal theoretic invariants}\label{Review_ideals}

The aim of this section is to briefly review some notions and tools from ideal theory. 

Let $(R,\mathfrak{m})$ denote a Notherian local ring and its maximal ideal (respectively, a standard graded ring over a field and its irrelevant ideal).
For an ideal $I\subset \mathfrak{m}$ (respectively, a homogeneous ideal $I\subset \mathfrak{m}$), the \emph{special fiber} of $I$ is the ring $\mathcal{F}(I):=\mathcal{R}(I)/\mathfrak{m}\mathcal{R}(I)$.
Note that this is an algebra over the residue field of $R$.
The (Krull) dimension of this algebra is called the \emph{analytic spread} of $I$ and is denoted $\ell (I)$. 

Quite generally, given ideals $J\subset I$ be ideals in a ring $R$,  $J$ is said to be a \emph{reduction} of $I$ if there exists an integer $n\geq 0$ such that $I^{n+1}=JI^n.$
An ideal shares the same radical with all its reductions.
Therefore, they share the same set of minimal primes and have the same codimension.
Obviously, any ideal is a reduction of itself, but one is interested in ``minimal'' reductions.
A reduction $J$ of $I$ is called \emph{minimal} if no ideal strictly contained in $J$ is a reduction of $I$.
The \emph{reduction number} of $I$ with respect to a reduction $J$ is the minimum integer $n$ such that $JI^{n}=I^{n+1}$. It is denoted by $\mathrm{red}_{J}(I)$. The (absolute) \emph{reduction number} of $I$ is defined as $\mathrm{red}(I)=\mathrm{min}\{\mathrm{red}_{J}(I)~|~J\subset I~\mathrm{is}~\mathrm{a}~\mathrm{minimal}~\mathrm{reduction}~\mathrm{of}~I\}.$ 
If $R/\mathfrak{m}$ is infinite, then every minimal reduction of $I$ is minimally generated by exactly $\ell(I)$ elements. In particular, every reduction of $I$ contains a reduction generated by $\ell(I)$ elements.

The following invariants are related in the case of $(R,\mathfrak{m})$:
$$\mathrm{ht}(I)\leq \ell(I) \leq \min\{\mu(I), \mathrm{dim}(R)\},$$
where $\mu(I)$ stands for the minimal number of generators of $I$.
If the rightmost inequality turns out to be an equality, one says that $I$ has maximal analytic spread.
By and large, the ideals considered in this work will have $\dim R\leq \mu(I)$, hence being of maximal analytic spread means in this case that $\ell(I)=\dim R$.

One next turns to notions that are more dependent on the choice of systems of generators.
Thus, quite generally, if $I\subset R$ is an arbitrary ideal in a Noetherian ring, picking a system of generators of $I$ gives rise to a free presentation 
$$R^m\stackrel{\phi}\lar R^n\rar I\rar 0,$$
where by abuse $\phi$ also denotes a matrix representative of the $R$-module map $R^m\rar R^n$.
Applying the symmetric functor it obtains an algebra presentation of $\mathcal S_R(I)$
$$\mathcal S_R(R^m)\stackrel{\mathcal S_R(\phi)}\lar \mathcal S_R(R^n)\rar \mathcal S_R(I)\rar 0,$$
over the polynomial ring $\mathcal S_R(R^n)\simeq R[\tt]:=R[t_1,\ldots,t_n]$, where $R^n=\sum_{i=1}^n Rt_i$.
In this passage, the image of $\mathcal S_R(\phi)$ in $R[\tt]$ is $I_1(\tt\cdot\phi)$, the ideal generated by the entries of the matrix product $\tt\cdot\phi$.
We call the latter a {\em presentation ideal} of the symmetric algebra of $I$.

One introduces in a similar way a {\em presentation ideal} of the Rees algebra $\mathcal R_R(I)=R[IT]\subset R[T]$, namely, as the kernel $\mathcal J$ of the surjective $R$-algebra homomorphism
$$R[\tt]\surjects R[IT], t_i\mapsto a_iT,$$
where $\{a_1,\ldots,a_n\}$ is the chosen system of generators.
It is fairly obvious that the latter homomorphism factor through the two surjective homomorphisms $R[\tt]\surjects \mathcal S_R(I)\surjects \mathcal R_R(I)$, hence the inclusion $I_1(\tt\cdot\phi)\subset \mathcal J$.
If further $I$ admits a regular element -- which will always be the case in this work -- then the Rees algebra $\mathcal R_R(I)$ is obtainable moding out the $R$-torsion ideal of $\mathcal S_R(I)$.
In particular, in this situation, $\mathcal J=I_1(\tt\cdot\phi):I^{\infty}$.

In the present ambient, one can add a further approximation to a system of generators of the presentation ideal $\mathcal J$ of the Rees algebra of $I$.
Namely, it is well-known or easy to see that $\mathcal J$ contains as subideal a presentation ideal $Q$ of the special fiber algebra $\mathcal R_R(I)/\fm \mathcal R_R(I)$.
Therefore, $(I_1(\tt\cdot\phi), Q)\subset \mathcal J$ and one says that $I$ (or its Rees algebra) is of {\em fiber type} if this inclusion is an  equality.

\section{Basic notions}

\subsection{The chaos invariant}

For convenience we introduce the following terminology.

\begin{Definition}\label{chaos_defined}\rm
	Let  $\phi$ denote an $n\times (n-1)$ matrix such that its entries are linear forms in the polynomial ring $R:=k[x,y,z].$ 
	Assume that $\hht I_1(\phi)=3$ and $\hht I_{n-1}(\phi)=2$.
	Then there is a uniquely defined integer $1\leq u\leq n-2$ such that
	\begin{equation}
	\hht I_t(\phi)=\left\{\begin{array}{cc}
	3,&1\leq t\leq u\\
	2,&u+1\leq t\leq n-1
	\end{array}\right.
	\end{equation} 
	We call $u=u(\phi)$ the {\em chaos invariant} of $\phi$.	
\end{Definition}
One can introduce a local version of this invariant in the following way: we know that $I_{n-1}(\phi)$ is a perfect ideal, hence its associated primes have codimension $2$. Further, for any $1\leq t\leq n-2$ one has $I_{n-1}(\phi)\subset I_t(\phi)$.  Let $p$ stand for any associated prime of $I_{n-1}(\phi)$. Then there exists a unique integer $1\leq u_{p}\leq n-2$ such that $p\in {\rm Min}(I_{u_{p}+1})\setminus {\rm Min}(I_{u_{p}})$.
We will call $u_p$ the {\em chaos invariant of $\phi$ at} $p$.

Clearly, one has
\begin{equation}\label{chaos_as_minimum}
u(\phi)=\min\{u_{p}\,|\,p\in{\rm Ass}(R/I_{n-1}(\varphi))\}
\end{equation} 

Recall that the group ${\rm Gl}_3(k)\times {\rm Gl}_n(k)\times {\rm Gl}_{n-1}(k)$ acts on the set of linear $n\times (n-1)$ matrices over $k$ through change of variables in $R$, elementary row operations and elementary column operations, respectively.
Two such matrices will be said to be {\em conjugate} if they belong to the same orbit of this action (see \cite[Lecture 9, Linear Determinantal Varieties in General]{Harrisbook}).

It is clear that the chaos invariants introduced above do not change under conjugation.

In the following we will systematically set $I:=I_{n-1}(\phi)$.
The basic result about the chaos invariants goes as follows.

\begin{Lemma}\label{chaos_basic}
Let  $\phi$ denote an $n\times (n-1)$ matrix whose entries are linear forms in the polynomial ring $R:=k[x,y,z]$, such that $I_1(\phi)=(x,y,z)$ and $\hht I=2$. Let $p$ stand for an associated prime of $I$.
Then:
\begin{enumerate}
	\item [{\rm (a)}] The rank of $\phi$ over $R/p$ is $u_p$.
	\item[{\rm(b)}] $\mu (I_{p})=n-u_{p}.$ 
	\item[{\rm(c)}] $I_{p}\subset {p_{p}}^{n-u_{p}-1}.$
\end{enumerate}
\end{Lemma}
\demo  (a) By definition of $u_p$, $p$ is not a minimal prime of $I_{u_{p}}$, hence $I_{u_{p}}\not\subset p$. Therefore, there exists an $u_p$-minor of $\phi$ not contained in $p$, i.e., the rank of $\phi$ over $R/p$ is at least $u_p$. On the other hand, still by definition,  $I_{u_p+1}(\varphi)\subset p$. Therefore, $I_t(\varphi)\subset p$ for every $t\geq u_p+1$. Thus, the rank of $\phi$ over $R/p$ is at most $u_p$.

(b) One has:
\begin{eqnarray*}\nonumber
\mu(I_{p})&=&\mu((I/pI)_p) \;\; ({\rm Nakayama})\\ \nonumber
&=&n-{\rm rk}_{R_{p}/{p_p}}(\varphi)=n-{\rm rk}_{R/p}(\varphi)\\ \nonumber
&=& n-u_p \;\; (\text{\rm by (a)})
\end{eqnarray*}

(c) This is an immediate consequence of knowing that $I_{u_{p}}\not\subset p$.
For then, locally at $p$, some $u_p$-minor of $\phi$ is invertible, so that up to conjugation $\phi$ has the form
\begin{equation}\label{block_diagonal}
\left(
\begin{array}{c@{\quad\vrule\quad}c}
\raise5pt\hbox{$\mathbb{I}$}&\raise5pt\hbox{$0$}\\[-6pt]
\multispan2\hrulefill\\
0&\psi
\end{array}
\right)
\end{equation}
where $\mathbb{I}$ is the identity matrix of order $u_p$ and $\psi$ is an $(n-u_p)\times (n-u_p-1)$ matrix with entries in $p_p$.
Therefore,  any generator of $I_p$ belongs to the power ${p_{p}}^{n-u_{p}-1}$.
\qed

\begin{Remark}\rm
	Although $I_p$ above has the same number of minimal generators as ${p_{p}}^{n-u_{p}-1}$, it may be the case that they are not all strictly of order $n-u_p-1$ -- see Theorem~\ref{Lan1} (b) below for the case where $u_p=1$.  
\end{Remark}

Our next concern has to do with the value of the analytic spread $\ell(I)=\dim \mathcal{F}(I)=\dim k[I_{n-1}]$.
Clearly, $2\leq \ell(I)\leq \dim R=3$.
Since $I$ is in general not generically a complete intersection, Cowsik--Nori theorem is not applicable to deduce that $\ell(I)\neq 2$, and neither is the previous work of Ulrich--Vasconcelos and Morey--Ulrich (\cite[Corollary 4.3]{UV}, \cite[Theorem 1.3]{MU}) because requiring that $I$ satisfy the condition called $G_3$ is tantamount to requiring that $I$ be generically a complete intersection.

Instead, by a quirk using a more encompassing version of a birationality criterion, we can show that $\ell(I)=3$ and more:

\begin{Theorem}\label{Fitting_main}
Let  $\phi$ denote an $n\times (n-1)$ matrix such that its entries are linear forms in the polynomial ring $R:=k[x,y,z].$ 
Assume that $\hht I_1(\phi)=3$ and $\hht I_{n-1}(\phi)=2$.
Set $I=I_{n-1}(\phi)$.
Then:
\begin{enumerate}
	\item[{\rm (a)}] The rational map $\mathfrak{F}$ defined by the maximal minors of $\phi$ is a birational map of $\pp^2$ onto its image.
	In particular, the ideal $I$ has maximal analytic spread.
	\item[{\rm (b)}] The inverse map to  $\mathfrak{F}$ in {\rm (a)} is defined by forms of degree $2$.
	\item[{\rm (c)}] ${\rm depth}(R/I^2)=0$.
\end{enumerate}
\end{Theorem}
\demo Set $u:=u(\phi)$ and let $p$ denote an associated prime of $I_{u+1}(\phi)$.
By Lemma~\ref{chaos_basic} (a), there is an $u$-minor of $\phi$ not contained in $I_u(\phi)$.
Therefore, up to conjugation we may assume that $p:=(y,z)$  and that $\phi$ has the following shape
\begin{equation}\label{phi_under_conjugation}
\phi=\left(\begin{array}{cccc|ccccc}
x+a_{1}&&&&&&\\
&x+a_{2}&&&&&\\
&&\ddots&&&&&\\
&&&x+a_{u}&&&\\
\hline
&&&\\
&&&\\
&&&
\end{array}\right),
\end{equation}
where $a_1,a_2,\ldots,a_{u}$ and the blank entries are linear forms in $y,z$ solely.

(a) In order to prove the statement, we use the criterion of \cite[Theorem 2.18 (b)]{AHA} in terms of a certain Jacobian dual matrix.
Namely, consider new variables $\tt=\{t_1,\ldots,t_n\}$ over $R$ and
let $B$ denote the uniquely defined matrix with entries in $k[\tt]=k[t_1,\ldots,t_n]$ satisfying the equality 
\begin{equation}\label{jacobian_duality}
(\tt)\cdot\phi=(x\, y\, z)\cdot B.
\end{equation}
Clearly, $B$ turns out to be a $3\times (n-1)$ matrix with linear entries. For convenience we look at its transpose: 
$$B^t=\left(\begin{array}{ccc}
t_1&\ell_{1,2}&\ell_{1,3}\\
\vdots&\vdots&\vdots\\
t_u&\ell_{u,2}&\ell_{u,3}\\
0&\ell_{u+1,2}&\ell_{u+1,3}\\
\vdots&\vdots&\vdots\\
0&\ell_{n-1,2}&\ell_{n-1,3}
\end{array}\right),$$
which is a submatrix of the so-called Jacobian dual matrix of $I$ as introduced in \cite[Section 2.3]{AHA}.
Note that the inequality $u\leq n-2$ implies the existence of at least one zero along the first column of $B^t$.

On the other hand $\ell_{u+1,2}\neq 0$.
This is because otherwise, $\phi$ would admit an entire column with entries depending only on $z$.
Since the maximal minors  belong to the ideal $(z)$ generated by the entries of that column this would contradict the standing assumption that $\hht I_{n-1}(\phi)=2$.

On the other hand, this same assumption implies that the maximal minors are $k$-linearly independent, hence admit no polynomial relation of degree $1$ in $k[\tt]$.
But the ideal $P$ of these relations is prime because it is the homogeneous defining ideal of the image of the rational map defined by the maximal minors.
This finally gives
 $$\det\left(\begin{array}{cc}t_u&\ell_{u,2}\\0&\ell_{u+1,2}\end{array}\right)\not\equiv 0\mod P,$$
telling us that the matrix $B^t$ has rank at least $2$ over $R[\tt]/P$. Therefore, so does the more encompassing Jacobian dual matrix as defined in \cite[Section 2.3]{AHA}. Moreover, by \cite[Corollary 2.16]{AHA} this rank is at most $\dim R-1=2$. Therefore, the Jacobian dual matrix has rank exactly $2$.
This is condition (b) in the aforementioned criterion.

(b) We draw on \cite[Theorem 2.18, Supplement (ii)]{AHA}, by which a representative of the inverse map can be taken to be the $2\times 2$ minors of an arbitrary $2\times 3$ submatrix of rank $2$ of the Jacobian dual matrix. In particular, one can choose such a matrix to be the horizontal slice
$$\left(\begin{array}{ccc}
t_u&\ell_{u,2}&\ell_{u,3}\\
0&\ell_{u+1,2}&\ell_{u+1,3}
\end{array}\right)$$
of $B^t$.

(c) Since the inverse map to $\mathfrak{F}$ is defined by forms of degree $2$ as in (b), then \cite[Proposition 1.4 (a)]{SymBir} tells us that there is an element in $I^{(2)}\setminus I^2$.
Clearly, then ${\rm depth}(R/I^2)=0$ since $I^{(2)}$ is the unmixed part of the second power.
\qed

\medskip

\begin{Proposition}\label{chaos_and_minimal_primes}
Let  $\phi$ denote an $n\times (n-1)$ matrix such that its entries are linear forms in the polynomial ring $R:=k[x,y,z]$, with $\hht I_1(\phi)=3$ and $\hht I_{n-1}(\phi)=2$. 

Setting $u=u(\phi)$, one has:
\begin{enumerate}
	\item[{\rm(a)}] If $n\geq 2(u+1)$, all the ideals $I_t(\phi)$ for $u+1\leq t\leq n-(u+1)$ have one single and the same minimal prime.
\end{enumerate}	
	Let $q$ denote the uniquely universal minimal prime in {\rm (a)}. Then:
\begin{enumerate}
	\item[{\rm(b)}]  For every associated prime of $I$ other than $q$ one has $u_p\geq n-(u+1)$ and $\mu(I_p)\leq u+1$.
\end{enumerate}
\end{Proposition}
\demo (a) Under the assumption that $n\geq 2(u+1)$, consider the chain of prime ideals 
$$I\subset I_{n-(u+1)}(\phi)\subset I_{n-(u+2)}(\phi)\subset \cdots\subset \underbrace{I_{u+1}(\phi)}_2\subset \underbrace{I_u(\phi)}_3,$$
where the tagged numbers denote codimensions arising from the definition of $u=u(\phi)$.
Fix a minimal prime $q$ of $I_{u+1}(\phi)$. Then $q$ is both a minimal prime of $I_{n-(u+1)}(\phi)$ and of $I$. By the definition of $u_q$, we have $u=u_q$.
Then, by Lemma~\ref{chaos_basic} (a) the rank of $\phi$ over $R/q$ is $u_q=u$. Therefore, up to conjugation $\phi$ has the form of (\ref{phi_under_conjugation}),
where $q=(y,z)$, and $a_1,\ldots,a_u$ and the blank entries are linear forms involving only $y,z$.

Let $J\subset I_{n-(u+1)}(\phi)$ denote the ideal generated by the $(n-(u+1))$-minors of the rightmost $n-(u+1)$ columns. Clearly, $q=(y,z)$ is the unique minimal prime of $J$ and $I\subset J$. Therefore, $J$ has height $2$ and hence $q$ is the unique minimal prime of $I_{n-(u+1)}(\phi)$ as well.

By tracing up the above chain of ideals $I_t(\phi)$, it is clear that $q$ is also the unique minimal prime of each of them for $u+1\leq t\leq n-(u+1)$.

(b) By assumption and item (a), $p$ is a minimal prime of $I$ but not of $I_{n-(u+1)}$.
Therefore, $u_p\geq n-(u+1)$ by definition.
Then, by Lemma~\ref{chaos_basic} (b) one has $\mu(I_p)\leq u+1$. 
\qed

\begin{Remark}\rm
The extreme case where $u(\phi)=1$ has been treated before in \cite{Lan} -- here all the ideals of minors, except $I_1(\phi)$ and $I$, have a unique and the same minimal prime. A full consideration of the results obtained in that paper will be taken up in a subsequent section.
\end{Remark}

\subsection{The special fiber and the reduction number}

Let $I\subset R:=k[\xx]$ stand for a homogeneous ideal in a standard graded ring over a field.
A very fundamental algebra related to  $I$ is its special fiber (algebra) $\mathcal{F}(I)$, as explained in Subsection~\ref{Review_ideals}.

The following result was observed in \cite[Proposition 1.2]{blup}:
\begin{Proposition}\label{CMfiber}
Let $I\subset R$ denote a homogeneous ideal such that its special fiber $\mathcal{F}(I)$ is Cohen--Macaulay. Then the reduction number $r(I)$ of $I$ coincides with the Castelnuovo--Mumford regularity ${\rm reg}(\mathcal F(I))$ of $\mathcal{F}(I)$.
\end{Proposition}
Here we think of  $\mathcal F(I)$ as a graded $S$-module, where $S=k[y_1,\ldots,y_n]$ with $n=\mu(I)$.

The proof rests on three ingredients: first, \cite[Proposition 1.85]{Wolmbook3} which tells us that when the special fiber is Cohen-Macaulay, one can read $r(I)$ off the Hilbert series of  $\mathcal F(I)$ as the degree of the polynomial in the numerator in its fractional form (the so-called $h$-{\em polynomial}).
Second, again since $\mathcal F(I)$ is Cohen--Macaulay, it has a minimal graded $S$-free resolution of length  $m:=n-\ell(I)$, and in addition the regularity is attained at the tail of the resolution, i.e., ${\rm reg}(\mathcal F(I))=\alpha-m$, where $\alpha$ is the largest shift in the minimal graded free resolution.
Finally,  the additivity of the Hilbert series along short exact sequences and the equality $\displaystyle HS(S^u(-v),t)=\frac{t^v}{(1-t)^n}\,u$, together give $b+m=\alpha$, where $b$ denotes the degree of the $h$-polynomial.
\qed

\medskip
Next is a remake of \cite[Proposition 3.8]{fat2} in greater generality.

\begin{Corollary}\label{red_at_most_2} Let $I\subset R$ denote a homogeneous ideal such that its special fiber $\mathcal{F}(I)$ is Cohen-Macaulay. Setting $n:=\mu(I)$, the following conditions are equivalent:
\begin{enumerate}
	\item[{\rm (i)}] The reduction number $r(I)$ is at most $\ell(I)-1$.
	\item[{\rm (ii)}] The largest shift in the minimal graded resolution of $\mathcal{F}(I)$ over $S=k[y_1,\ldots,y_n]$ is at most $n-1$.
	\item[{\rm (iii)}] The Hilbert function $H_{\mathcal{F}(I)}(t)$ of $\mathcal{F}(I)$ coincides with its Hilbert polynomial for all $t\geq 1$.
\end{enumerate}
\end{Corollary}
\demo
The equivalence of (i) and (ii) follows immediately from Proposition~\ref{CMfiber} since the homological dimension of $\mathcal{F}(I)$ is $n-\ell(I)$. 
On the other hand, again by \cite[Proposition 1.85]{Wolmbook3}, $r(I)$ is the degree of the $h$-polynomial in the Hilbert series of $\mathcal{F}(I)$.
Therefore, since $\dim \mathcal{F}(I)=\ell(I)$, the equivalence of (i) and (iii) follows from \cite[Proposition 4.3.5 (c)]{BH}.
\qed

\begin{Remark}\label{red_one}\rm
	In the above corollary, $r(I)$ is frequently said to have the expected value if equality takes place in item (i). 
	The extreme lowest value $r(I)=1$ corresponds of course to the event that $\mathcal{F}(I)$ has minimal multiplicity (degree), hence the corresponding projective variety is either a rational normal scroll or a cone over a Veronese.
\end{Remark}

When $R=k[x,y,z]$, an additional restriction on the reduction numbers of the localizations will give that $r(I)\leq 2$ is in fact equivalent to having $\mathcal{F}(I)$  Cohen--Macaulay.


\begin{Theorem}\label{red_no2} Let $I\subset R=k[x,y,z]$ denote a codimension $2$ perfect homogeneous ideal satisfying the following conditions:
\begin{enumerate}
	\item[{\rm (a)}] $I$ is linearly presented
	\item[{\rm (b)}] $r(I_p)\leq 1$ for every associated prime $p$ of $R/I$
\end{enumerate}
	Then the following conditions are equivalent:
	\begin{enumerate}
		\item[{\rm(i)}] $r(I)\leq 2.$
		\item[{\rm(ii)}] The Rees algebra $\mathcal{R}_R(I)$ is Cohen--Macaulay.
		\item[{\rm(iii)}] The special fiber $\mathcal{F}(I)$ is Cohen--Macaulay.
	\end{enumerate}
\end{Theorem} 
\demo
(i) $\Rightarrow$ (ii)
We draw on the following result of Cortadellas--Zarzuela, as slightly reformulated in \cite[Theorem 3.6]{fat2}:
\begin{Proposition}\label{CortZar}
	Let $R$ denote a Cohen--Macaulay local ring and let $I\subset R$ stand for an ideal satisfying the following properties:
	\begin{enumerate}
		\item[{\rm (1)}]  $I$ is unmixed, and $R/I$ and $R/I^2$ have different depths.
		\item[{\rm (2)}]  $\max\{r(I_P)\,|\, J\subset P, \hht(P)=\hht(I)\}\leq 1$.
		\item[{\rm (3)}] $\ell(I)=\hht(I)+1\geq 3$.
		\item[{\rm (4)}] $r(I)\leq 2$.
	\end{enumerate}
	Then ${\rm depth}\, \mathcal{R}_R(I)=\min\{{\rm depth} \,R/I, {\rm depth} \,R/I^2+1\} +\hht(I)+1$.
\end{Proposition}
To apply this result, thus concluding that ${\rm depth}\, \mathcal{R}_R(I)=4=\dim R+1$, we note that conditions (1) and (3) follow from the standing hypotheses (specially, (a)) by drawing upon Theorem~\ref{Fitting_main}, while condition (2) is contained in hypothesis (b).

(ii) $\Rightarrow$ (iii).
Since $I$ is generated in a fixed degree, then $\mathcal{F}(I)$ is isomorphic as graded $k$-algebra to the $k$-subalgebra $k[I_{n-1}t]\subset R[It]$, where $n=\mu(I)$. Therefore, the result follows from \cite[Proposition 3.10]{fat2}.

(iii) $\Rightarrow$ (i) We apply Corollary~\ref{red_at_most_2} by showing its condition (iii) holds.
To see this, we first show that this condition holds locally at every associated prime $p$ of $R/I$. But locally at such $p$ the ideal $I_p$ is $p_p$-primary in $R_p$. Since $r(I_p)\leq 1$ by (b), \cite[Theorem 2.1]{Huneke} implies that the Hilbert--Samuel function $HS_{I_p}(t)$ of $I_p$ is polynomial for any $t\geq 1$.

Let $e(M)$ denote the multiplicity (degree) of a graded module $M$ over $R$.
Letting $t\geq 1$, the associativity formula for multiplicities yields
$$e(R/I^t)=\sum_{p\in{\rm Min}(R/I)} \lambda(R_{p}/I_{p}^t)e(R/p)=\sum_{p\in{\rm Min}(R/I)} \lambda(R_{p}/I_{p}^t)=\sum_{p\in{\rm Min}(R/I)}HS_{I_p}(t).$$
It follows that $e(R/I^t)$ has the values of a polynomial for every $t\geq 1$.

But the graded  Hilbert function of $\mathcal{F}(I)$ has the form
$$H_{\mathcal{F}(I)}(t)=\mu(I^t)={\mu(I)\,t+1\choose 2}-e(R/I^t),$$ 
for every $t\geq 1$, where $\mu(\_)$ denotes minimal number of generators.
Therefore, the Hilbert function of $\mathcal{F}(I)$ has the values of a polynomial in every degree $\geq 1$, hence coincides with its Hilbert polynomial throughout.
\qed

\section{The role of a Jacobian dual submatrix}

\subsection{A $1$-generic submatrix of the Jacobian dual matrix}

The matrix $B$ obtained in (\ref{jacobian_duality}) is useful not only as a major vehicle to test for birationality and to obtain a representative of the inverse map as described in \cite{AHA}, but often to guess generators of the homogeneous defining ideal of the image of the map.  

For convenience, we write it explicitly:

\begin{equation}\label{jacobian_dual_transpose}
B=\left(\begin{array}{ccccccc}
t_1&t_2&\ldots&t_u&0&\ldots&0\\
\ell_{1,2}&\ell_{2,2}&\ldots&\ell_{u,2}&\ell_{u+1,2}&\ldots&\ell_{n-1,2}\\
\ell_{1,3}&\ell_{2,3}&\ldots&\ell_{u,3}&\ell_{u+1,3}&\ldots&\ell_{n-1,3}
\end{array}\right),
\end{equation}
where $\ell_{i,2}, \ell_{i,3}$ are certain linear forms in $k[\tt]=k[t_1,\ldots,t_n]$.

Set
\begin{equation}\label{2x2piece}
B':=\left(\begin{array}{ccccccc}
\ell_{u+1,2}&\ldots&\ell_{n-1,2}\\
\ell_{u+1,3}&\ldots&\ell_{n-1,3}
\end{array}\right)
\end{equation}
We have already remarked in the proof of  Theorem~\ref{Fitting_main} (b) that no entry of this matrix vanishes.
We prove more:

\begin{Proposition}\label{genericity} The matrix $B'$ is $1$-generic.
\end{Proposition}
\demo We are to show that $B'$ acquires no zero entry under the action of ${\rm GL}_2(k)\times {\rm GL}_{n-1-u}(k)$.
Let us argue by contradiction, assuming that there exist
$A\in {\rm GL}_2(k)$ and $C\in {\rm GL}_{n-1-u}(k)$ such that $A\cdot B'\cdot C$ has a null entry. 

Set 
$$\tilde{A}:=\left(\begin{array}{cc}1&\\&A\end{array}\right)\in{\rm GL}_{3}(k)\quad\quad\mbox{and}
\quad\quad\tilde{C}=\left(\begin{array}{cc}\mathbb{I}_u&\\
&C\end{array}\right)\in {\rm GL}_{n-1}(k),$$
where $\mathbb{I}_u$ denotes the identity matrix of order $u$.

Then the transform of $B$ by the latter matrices has the shape
\begin{equation}\label{transform_of_B}
\tilde{A} B\tilde{C}=
\left(\begin{array}{ccc|cccccc}
t_1&\ldots&t_u&\boldsymbol0\\
\ast&\ldots&\ast&A B' C
\end{array}\right),
\end{equation}
where the rightmost block has a columns with at least two zero entries.

{\sc Claim:} (\ref{transform_of_B}) returns the Jacobian dual matrix on the right side of (\ref{jacobian_duality}) relative to the matrix $\tilde{\phi}$ obtained by acting on $\phi$ with matrices $\tilde{A}$ and $\tilde{C}$.

The assertion is clear regarding the operation by $\tilde{C}$ since one just have to multiply both members of (\ref{jacobian_duality}) on the right by $\tilde{C}$.

The argument regarding the left-side operation is more delicate since for $\phi$ one means change of variables in $R$, while for $B$ one has to engage on a direct verification.

Rewriting $\phi$ in new coordinates $x',y',z'$ given by
$$\left(\begin{array}{c}
x'\\
y'\\
z'
\end{array}\right)=(\tilde{A}^{-1})^t\left(\begin{array}{c}
x\\
y\\
z
\end{array}\right),$$
yields
$$\tt\cdot\varphi\tilde{C}=(x'\,y'\,z')\cdot\tilde{A} B \tilde{C}.$$

This way, if $A\cdot B'\cdot C$ has a zero entry then there exists an entire column of  $\tilde{A} B \tilde{C}$ with two zeros. Therefore, $\varphi\tilde{C}$ has a column depending solely in one of the two variables $y'$ ou $z',$ which would imply that   $I_{n-1}(\varphi\tilde{C})$ is contained in the principal ideal generated by one of these variables  -- an absurd as  $I_{n-1}(\varphi\tilde{C})$ has codimension $2$.  

As a consequence of the claim,  $\tilde{\phi}$ has an entire column depending on one single variable. As already explained, this contradicts that ${\rm cod}\, I_{n-1}(\tilde{\phi})={\rm cod}\, I_{n-1}(\phi)=2$.
\qed

\begin{Theorem}\label{quadrics_is_scroll}
Let $\phi$ stand for the matrix in {\rm (\ref{phi_under_conjugation})} and let $B'$ be the associated matrix as in {\rm (\ref{2x2piece})}. Then $I_2(B')$ is isomorphic to the homogeneous coordinate ring of a rational normal scroll of dimension $u+1$ in $\pp^{n-1}$, where $u$ denotes the chaos invariant of $\phi$. In particular, it is a Cohen--Macaulay normal domain.
\end{Theorem}
\demo This is a consequence of \cite[Proposition 9.12]{Harrisbook}, or rather from its subsequent generalization due to F. Schreyer,   where it is shown that any  $1$-generic $2\times r$ matrix, with $r\leq n-1$, is conjugate to a matrix of the shape
$$\left(\begin{array}{cccc|cccc|c|cccc}
x_0 & x_1 & \ldots & x_{a_1-1} & y_{0} & y_{1} & \ldots & y_{a_2-1} & \ldots & z_0& z_1 &\ldots & z_{a_s-1}\\
x_1 & x_2 & \ldots & x_{a_1} & y_{1} & y_{2} & \ldots & y_{a_2} & \ldots & z_1& z_2 &\ldots & z_{a_s}
\end{array}
\right) $$
formed by $s$ mutually independent Hankel blocks, with $r=a_1+\cdots +a_s$.
The $2$-minors generate the homogeneous defining ideal of a rational normal scroll in projective space $\pp^{r+s-1}$, hence it is a prime ideal and the corresponding homogeneous coordinate ring is normal.

We also note that the rational normal scroll has dimension $s$.
Applying to $B'$, with $r=n-u-1$, one has that $I_2(B')$ has the expected codimension $n-u-1-2+1=n-u-2$ in $k[\tt]=k[t_1,\ldots,t_n]$, hence $\dim  k[\tt]/I_2(B')=u+2$ and therefore $s=u+1$ in our case.
Then $k[\tt]/I_2(B')$ is isomorphic to the homogeneous coordinate ring of a rational normal scroll of dimension $u+1$ in $\pp^{n-1}$.
\qed

\subsection{The fiber type property}

\subsubsection{Arbitrary chaos invariant}

In this part we turn to the Rees algebra $R[IT]\subset R[T]$ of the ideal $I=I_{n-1}(\phi)$.
We refer to the same variables $\tt=\{t_1,\ldots,t_n\}$ introduced in the previous subsection and to the relation (\ref{jacobian_duality}).

Let $\mathcal{J}$ denote the presentation ideal of $R[IT]$ over $R[\tt]$, i.e., the kernel of the $R$-homomorphism $R[\tt]\rar R[IT]$ mapping $t_i$ to $\Delta_iT$, where $\Delta_i$ is the maximal minor of $\phi$ deleting the $i$th row.

One knows that $I_1((\tt)\cdot\phi)$ is the presentation ideal of the symmetric algebra of $I$, hence $I_1((\tt)\cdot\phi)\subset \mathcal{J}$.
From relation (\ref{jacobian_duality}), an easy consequence of the Cramer rule implies that $I_3(B)\subset \mathcal{J}$.
Therefore, one has 
$$(I_1((\tt)\cdot\phi), I_3(B))\subset (I_1((\tt)\cdot\phi),Q)\subset\mathcal{J},$$
where $Q$ denotes the presentation ideal of the special fiber of $I$ or, equivalently, the homogeneous defining ideal of the image of the birational map defined by the maximal minors $\Delta_i, 1\leq i\leq n$ of $\phi$.

Now looking at the shape of $B$ as in (\ref{jacobian_dual_transpose}), one sees that $I_3(B)=(I_3(B_1), \mathcal{B})$, with
\begin{equation}
B_1=\left(\begin{array}{cccc}
t_1&t_2&\ldots&t_u\\
\ell_{1,2}&\ell_{2,2}&\ldots&\ell_{u,2}\\
\ell_{1,3}&\ell_{2,3}&\ldots&\ell_{u,3}
\end{array}\right).
\end{equation}
and $\mathcal{B}$ is the sum of the ideals $I_2(B''_{i,j}),\, 1\leq i<j\leq u$, with
\begin{equation}
B''_{i,j}:=\left(\begin{array}{c|c}
D_{i,j} & B'
\end{array}\right),
\end{equation}
where $B'$ is as in (\ref{2x2piece}) and $D_{i,j}$ denotes the  column vector with entries the $2\times 2$ minors
\begin{equation}\label{2x2dets}
\det\left(\begin{array}{cc}t_i&t_j\\
\ell_{i,2}&\ell_{j,2}\end{array}\right) \quad {\rm and} \quad
\det\left(\begin{array}{cc}t_i&t_j\\
\ell_{i,3}&\ell_{j,3}\end{array}\right).
\end{equation}

\begin{Question}\rm Notation as above. An opening question is whether the ideal $ (I_3(B_1), \mathcal{B})$ has codimension $n-3$.
\end{Question}
The above ideal is a weak approximation to the ideal of the special fiber of $I$.

A neat case goes next.

\subsubsection{Chaos invariant $1$}

The following results contain the main findings in \cite{Lan} from a different perspective, and possibly in a simpler fashion.

We start with a preliminary assertion about matrices with binary linear entries.

\begin{Lemma}\label{jac_dual_binary} {\rm ($k$ algebraically closed)}
Given an integer $m\geq 1$, let $M$ denote an $m\times m$ matrix with linear entries in  $k[y,z].$ Let $C$ denote the uniquely defined $2\times m$ matrix whose entries are linear forms in the polynomial ring  $k[t_1,\ldots,t_m]$ satisfying the matrix equality
	$$(t_1 \cdots t_m)\cdot M=(y\,\,z)\cdot C.$$
Suppose that $\det M\neq 0$ and that $I_{m-1}(M)$ is $(y,z)$-primary. Then $I_2(C)$ has the expected codimension $m-2+1=m-1.$ 
\end{Lemma}
\demo
Firstly,  since $k$ is algebraically closed, $M$ is conjugate to a matrix of the shape
$$\widetilde{M}:=\left(\begin{array}{c|cc}
y&\boldsymbol{\alpha}\\
\hline
\mathbf 0&M'
\end{array}\right)$$
where $\det M'\neq 0$ and the entries of $\boldsymbol{\alpha}=(\alpha_2z\cdots \alpha_{n}z)$ depend only on $z$, $\alpha_i\in k$. 

This fact is most probably well-known being essentially a result about matrices in one single variable. For the sake of completeness, we give the following argument:
write $M=yA+zB$, for certain $m\times m$ matrices over $k$.
Since $\det M$ is a nonzero binary form over an algebraically closed field, it has a zero $(a,b)\in k^2$, where, say, $a\neq 0$. Substituting gives $Av=(b/a) Bv$, for suitable nonzero vector $v\in k^m$.
Letting $\{v_1=v,v_2,\ldots, v_m\}$ denote a vector basis of $k^m$, there exists an invertible $m\times m$ matrix $T$ over $k$ such that the first column of 
$TAT^{-1}$ is the first column of $TBT^{-1}$ times $(b/a)$. 
In this way we get that the entries of the first column of $TMT^{-1}$ are scalar multiples of $(y+(b/a)z)$.  Additional elementary row operations on $TMT^{-1}$  will make all entries of its first column vanish except one, as was to be shown.

Let us now switch to the proof of the main statement.
We induct on $m$.

Clearly, $M'$ is an $(m-1)\times (m-1)$ matrix satisfying the same hypothesis as $M$, namely, $\det M'\neq 0$ and $I_{m-2}(M)$ is $(y,z)$-primary. 
Writing, similarly
		$$(t_2 \cdots t_n)\cdot M'=(y\,\,z)\cdot C',$$
where $C'$ is a uniquely defined $2\times (m-1)$ matrix of linear forms in $k[t_2,\ldots,t_m]$, by the inductive hypothesis the ideal $I_2(C')$ has codimension $m-2$.
Say,
$$C':=\left(\begin{array}{ccccccc}
\lambda_{1,2}&\ldots&\lambda_{1,m}\\
\lambda_{2,2}&\ldots&\lambda_{2,m}
\end{array}\right),$$
with $\lambda_{ij}\in k[t_2,\ldots,t_m]$.

Then, taking in account the above shape of $\widetilde{M}$, $C$ turns out to be  conjugate to
\begin{eqnarray}\label{conjugate_sparse}
\widetilde{C}=\left(\begin{array}{ccccccc}
t_1&\lambda_{1,2}&\ldots&\lambda_{1,m}\\
0&\lambda_{2,2}+\alpha_2t_1&\ldots&\lambda_{2,m}+\alpha_{m}t_1
\end{array}\right)
\end{eqnarray} 
Therefore, it suffices to show that $I_2(\widetilde{C})$ has codimension $m-1$.

For this, first note that the linear forms $\lambda_{2,2}+\alpha_2t_1,\ldots,\lambda_{2n}+\alpha_mt_1$ are $k$-linearly independent.
Indeed, otherwise up to conjugation $M$ would have two columns involving only one of the two variables $y,z$ and that would imply that ${\rm cod} \, I_{m-1}(M)=1$.

Now, if $P\supset I_2(\widetilde{C})$ is any prime ideal then it must contain either the above linear forms, in which case its codimension is at leat $m-1$; or else it has to contain $t_1$, hence also $I_2(C')$.
Thus, $P$ contains $(I_2(C'), t_1)$ which has codimension $m-2+1=m-1$ by the inductive hypothesis and the fact that $t_1$ is nonzerodivisor thereof.
\qed

\begin{Proposition}\label{codim_minors_jac_dual}
Let  $\phi$ denote an $n\times (n-1)$ matrix such that its entries are linear forms in the polynomial ring $R:=k[x,y,z]$, with $\hht I_1(\phi)=3$ and $\hht I=2$.
Assume that $u(\phi)=1$ and let $B$ denote the matrix {\rm (\ref{jacobian_dual_transpose})} in this case.
Then the ideal $I_2(B)$ has codimension $n-1$.
\end{Proposition}
\demo Recall that
\begin{equation}\label{jacobian_dual_lan}
B=\left(\begin{array}{ccccccc}
t_1&0&\ldots&0\\
\ell_{2,1}&\ell_{2,2}&\ldots&\ell_{2,n-1}\\
\ell_{3,1}&\ell_{3,2}&\ldots&\ell_{3,n-1}
\end{array}\right),
\end{equation}
where $\ell_{2,i}, \ell_{3,i}$ are certain linear forms in $k[\tt]=k[t_1,\ldots,t_n]$.
From the shape of the matrix $\phi$ as in (\ref{phi_under_conjugation}), one can write
$$\varphi=\left(\begin{array}{cccccc}
L\\
\hline
M
\end{array}\right),$$
where $M$ is a $(n-1)\times (n-1)$ matrix with linear  entries in $k[y,z]$.
Let $D$ denote the unique matrix with entries in $k[t_2,\ldots,t_n]$ satisfying the equality
$$(t_2\cdots t_n)\cdot M=(y\; z)D.$$
In addition, $\det M\neq 0$ as it is one of the minimal generators of our main ideal $I=I_{n-1}(\phi)$. On the other hand, $I_{n-2}(M)$ is $(y,z)$-primary since it contains the ideal $I$ and hence has codimension $2$.
Therefore, one can apply Lemma~\ref{jac_dual_binary} with $m=n-1$ to conclude that $I_2(D)$ has codimension $n-2$, hence $(I_2(D), t_1)$ has codimension $n-1$.

On the other hand, letting $C$ denote the submatrix of $B$ of the last two rows, one has $(I_2(C), t_1)=(I_2(D), t_1)$ and hence the first also has codimension $n-1$.
But $I_2(B)=(I_2(C), t_1 \mathcal{L})$, where $\mathcal L$ denotes the set of entries of the of the submatrix
\begin{equation}\label{2x2piece_lan}
B':=\left(\begin{array}{ccccccc}
\ell_{2,2}&\ldots&\ell_{2,n-1}\\
\ell_{2,3}&\ldots&\ell_{3,n-1}
\end{array}\right)
\end{equation}
By Proposition~\ref{genericity} $B'$ is $1$-generic, hence $\mathcal L$ span a vector space of dimension at least $n-2+2-1=n-1$ (see \cite[Proposition 1.3]{Eisenbud2}). Therefore $(I_2(C),\mathcal L)=(\mathcal L)$ too has codimension at least $n-1$.

This proves the statement.
\qed

\medskip
We now get a grip on the main results for the case where the chaos invariant of $\phi$ is $1$ (minimal possible).

\begin{Theorem}\label{Lan1}
Let  $\phi$ denote an $n\times (n-1)$ matrix such that its entries are linear forms in the polynomial ring $R:=k[x,y,z]$, with $\hht I_1(\phi)=3$ and $\hht I=2$.
Assume that $u(\phi)=1$.
One has:
	\begin{enumerate}
		\item[{\rm (a)}] Let $q$ denote the unique minimal prime ideal of all the ideals $I_t(\phi)$, for $t\neq 1, n-1$, as obtained in {\rm Proposition~\ref{chaos_and_minimal_primes} (a)}. Then, for every associated prime $p$ of $I$ other than $q$, $I_p$ is a complete intersection.
		\item[{\rm(b)}] Let $q$ be as in {\rm (a)}. Then $I_{q}$ admits a set of $n-1$ generators $$\{g_1,\ldots,g_r,h_{r+1},\ldots,h_{n-1}\},\; \text{\rm for some}\;  1\leq r\leq n-1,$$ 
		 with $g_1,\ldots, g_r\in {q_q}^{n-2}\setminus {q_q}^{n-1}$ and $h_{r+1},\ldots,h_{n-1}\in{q_q}^{n-1}\setminus {q_q}^{n}.$
		 
		 In addition, $I_q={q_q}^{n-2}$ if and only if $r=n-1$.
		\item[{\rm (c)}] The special fiber $\mathcal{F}(I)$ is either the homogeneous coordinate ring of a rational normal scroll or of a cone over the Veronese embedding of $\pp^1$. 
		In particular, $\mathcal{F}(I)$  is Cohen--Macaulay, has minimal degree and $r(I)\leq 1$.
		\item[{\rm (d)}] The Rees algebra of $I$ is Cohen--Macaulay.
		\item[{\rm (e)}] The Rees algebra of $I$ is of fiber type and normal.
	\end{enumerate}
\end{Theorem}
\demo
(a) This follows immediately from Proposition~\ref{chaos_and_minimal_primes} (b).

(b) Up to conjugation, we may assume that $q=(y,z)$ and that $\phi$ is in the canonical form as in (\ref{phi_under_conjugation}):

$$\phi=\left(\begin{array}{c|c}
x+a&\\
\hline\\ [-6pt]
&\;H
\end{array}\right),$$
where $a$ and the blank entries are $1$-forms involving only $y,z$ and and so are the entries of the $(n-1)\times (n-2)$ matrix $H$. Since $x+a$ is invertible locally at $q$, the ideal $I_q$ is generated by the images of the maximal minors of $\phi$ fixing the first row. The latter have the form
\begin{equation}
(x+a)\Delta_i+D_i,
\end{equation}
for $1\leq i\leq n-1$, where $\Delta_i$ is the $(n-2)$-minor of $H$ excluding the $i$th row and $D_i\in q^{n-1}\setminus q^n.$
Since the $(n-1)$-minor of $\phi$ excluding the first row is nonzero and its Laplace expansion along the leftmost column is a linear combination of the $\Delta_i$'s, it follows that there is at least one nonvanishing $\Delta_i$. On the other hand, since $x+a$ is invertible in $R_q$ then $(x+a)\Delta_i+D_i\in {q_q}^{n-2}\setminus {q_q}^{n-1}$ if and only if $\Delta_i\neq 0$. 
This proves the required statement.

To get the additional statement, it suffices to argue that if $r=n-1$ then $I_q$ is generated by $n-1$ $k$-linearly independent elements of order $n-2$, hence the modules $R_q/I_q$ and $R_q/{q_q}^{n-2}$ have same length.

(c) Since $\mathcal{F}(I)$ has codimension $n-3$, it follows from Proposition~\ref{genericity} and Theorem~\ref{quadrics_is_scroll} that the homogeneous defining ideal of  $\mathcal{F}(I)$ is generated by the $2$-minors of a $1$-generic $2\times (n-2)$ matrix of linear forms.
By \cite[Proposition 9.4 and Proposition 9.12]{Harrisbook}, this matrix is conjugate to the defining matrix of a rational normal scroll or a cone over the standard rational normal curve in $\pp^{n-2}$. The latter are varieties of minimal degree, hence  $\mathcal{F}(I)$ is Cohen--Macaulay. The value $r(I)\leq 1$ follows from the formula of the degree of $\mathcal{F}(I)$ coming from its Hilbert series.

(d) We apply Theorem~\ref{red_no2}, using (b) above, in which case we only have to prove that $r(I_p)\leq 1$ for every associated prime of $i$. If $p\neq q$, this follows from (a) above.
We now deal with the case of the prime $q$.
Still from (a) we know quite generally that the length $\lambda(R_p/I_p^t)$ is a polynomial function for any $t\geq 1$.
On the other hand, by (c) we also know that the Hilbert function of $H_{\mathcal{F}(I)}(t)$ is polynomial for every $t\geq 1$.
Since
$$\lambda(R_q/I_q^t) = H_{\mathcal{F}(I)}(t) - {{nt+1}\choose {2}} +  
\sum_{p\in {\rm Ass}(R/I)\setminus \{q\} }\lambda(R_p/I_p^t),$$
it follows that $\lambda(R_q/I_q^t)$ is also polynomial for every $t\geq 1$.
Applying \cite[Theorem 2.1]{Huneke} we obtain that $r(I_q)\leq 1$.

(e) 
In the notation of the proof of Proposition~\ref{codim_minors_jac_dual}, we claim that $\mathcal I:=(I_1((\xx)\cdot B),I_2(B'))$ is the presentation ideal $\mathcal J$ of the Rees algebra of $I$ as introduced in the beginning of this subsection.

For this, it suffices to show that $\mathcal I$ is prime ideal of codimension $n-1$ and $\mathcal I\subset \mathcal J$.
First, the containment: clearly, $I_1((\xx)\cdot B)\subset \mathcal J$ as the first is the presentation ideal of the symmetric algebra of $I$. Next, we have seen that by an easy application of Cramer rule, one has $I_3(B)\subset \mathcal J$ as well.
But $I_3(B)=t_1 I_2(B')$ and $\mathcal J$ is a prime ideal not containing $t_1$.
Therefore, $I_2(B')\subset \mathcal J$ as well.

First recall from Theorem~\ref{quadrics_is_scroll} that $I_2(B')$ is the defining ideal of a rational normal scroll. On the other hand, one can write
$$ \mathcal I=(t_1x+\ell_{2,1}y+\ell_{3,1}z, I_2(B'')),
$$
where 
$$B'':=\left(\begin{array}{ccccccc}
z&\ell_{2,2}&\ldots&\ell_{2,n-1}\\
-y&\ell_{2,3}&\ldots&\ell_{3,n-1}
\end{array}\right)$$
Therefore, $I_2(B'')$ is also the defining ideal of a rational normal scroll. 
In particular, it is a Cohen--Macaulay prime ideal of codimension $n-2$. 
Now, $I_2(B'')$ being a prime generated in degree $2$ in the variables $y,z,\tt$ forces $t_1x+\ell_{2,1}y+\ell_{3,1}z$ to be a nonzerodivisor on it. It follows that $\mathcal I$ is Cohen--Macaulay of codimension $n-1$.
It remains to show that it is a prime ideal. We do this by showing that the ring $R[\tt]/\mathcal I$ is normal (and hence,  the Rees algebra of $I$ will be normal).

It suffices to show the condition $(R_1)$ of Serre's.
Note that the Jacobian matrix $\Theta$ of $\mathcal I$ has the following shape
$$\Theta=\left(\begin{array}{cccc}
B&0\\
\phi&\Theta'
\end{array}\right)$$
where $\Theta'$ denotes the Jacobian matrix of $I_2(B').$
We have to show that the Jacobian ideal $(I_{n-1}(\Theta), \mathcal I)$ has codimension at least $n-1+2=n+1$.
A calculation leads to the inclusion
$$(I_{n-3}(\Theta')\cdot I_{2}(B), I_2(B'),I)\subset (I_{n-1}(\Theta), \mathcal I).$$
Now, the ideal $(I_{n-3}(\Theta'), I_2(B'))$ has codimension at least $n+1$ since it is the Jacobian ideal of the normal ring $k[\tt]/I_2(B')$.
Therefore, it remains to prove that the ideal $(I_{2}(B), I_2(B'),I)=(I, I_2(B))$ has codimension at least $n+1$ or, equivalently, that $I_2(B)$ has codimension at least $n-1$.
But this is the statement of Proposition~\ref{codim_minors_jac_dual}.
\qed

\begin{Remark}\rm
(i) The shape of the generators in item (b) of the above theorem look sort of vague. However, any possibility with $1\leq r\leq n-1$ actually happens and in fact even when $I$ is a monomial ideal; for a typical example, one can take
$$I=(xy^{n-2},xy^{n-3}z,\ldots, xy^{n-(r+1)}z^{r-1}, y^{n-(r+1)}z^{r},y^{n-(r+2)}z^{r+1}\ldots, yz^{n-2}, z^{n-1}).$$
(See Section~\ref{monomial} for a more systematic way of dealing with the monomial case). 

(ii) As to the two alternatives in item (c), both can take place. For an emblematic example of the case where the special fiber is a cone over the Veronese embedding of $\pp^1$, consider the ideal $I=(x.(y,z)^{n-2}, z.f)$, where $f\in k[y,z]_{n-2}, f\neq 0$.
An easy calculation shows that the special fiber is the ring of a cone over the ring of the rational normal curve in $\pp^{n-2}$.
\end{Remark}

\section{Paradigm classes}

In this part we introduce three classes of linearly presented codimension $2$ ideals along with the chaos invariants of the respective syzygy matrices.

\subsection{Linearly presented ideals of plane fat points}

We assume that $k$ is an algebraically closed field.

An ideal of plane fat points has the form
$$I=p_1^{m_1}\cap\cdots\cap p_r^{m_r},$$
where $p_i\subset R=k[x,y,z]$ is a codimension $2$ prime ideal (hence, generated by two independent linear forms) and $m_i$ is a positive integer (called {\em multiplicity} of the point associated to the prime ideal $p_i$).
We assume that $m_1\geq\cdots\geq m_r$.

See \cite{fat2} for a class of such ideals that are linearly presented up to a quadratic transformation (i.e., a transformation induced by a plane quadratic Cremona map).

The following establishes the value of the chaos invariant for these ideals.
\begin{Proposition}\label{chaos_4_fat}
	Let $I=p_1^{m_1}\cap\cdots\cap p_r^{m_r}$ denote an ideal of fat points as above.
	Suppose that $I$ is linearly presented, with syzygy matrix $\phi$ of size $n\times (n-1)$.
	Then $u(\phi)=n-m_1-1$, where $m_1$ is the highest multiplicity.
\end{Proposition}
\demo
On one hand, for every $i$, one clearly has $\mu(I_{p_i})=\mu({p_i}_{p_i}^{m_i})=m_i+1$.

On the other hand, we know from  Lemma~\ref{chaos_basic} (b) that $\mu(I_{p_i})=n-u_{p_i}$, where $u_{p_i}$ denotes the chaos invariant of $\phi$ at $p_i$.
This gives $u_{p_i}=n-m_i-1$, for every $i$.
By (\ref{chaos_as_minimum}), $u(\phi)=\min\{n-m_i-1\,|\, 1\leq i\leq r\}=n-m_1-1$, as required.
\qed

\begin{Corollary}  Let $I=p_1^{m_1}\cap\cdots\cap p_r^{m_r}$ denote an ideal of fat points as above.
	Suppose that $I$ is linearly presented, with syzygy matrix $\phi$ of size $n\times (n-1)$.
	Then:
	\begin{enumerate}
		\item[{\rm a}] The minimal primes of $R/I_{u+1}$ are the primes of $I$ having the highest multiplicity.
		\item[{\rm b}] The following conditions are equivalent:
		\begin{itemize}
			\item[{\rm (i)}] $u(\phi)=n-2$
			\item[{\rm (ii)}] $I$ has only simple points {\rm (}i.e.,  $I$ is a radical ideal{\rm )} 
			\item[{\rm (iii)}] $I$ is generically a complete intersection
			\item[{\rm (iv)}] $I$  satisfies property $G_3$.
		\end{itemize}
	\end{enumerate}
\end{Corollary}
\demo (a) This follows from the proof of the previous proposition. 

(b) By the previous proposition, (i) and (ii) are obviously equivalent.
Since a radical ideal is generically a complete intersection, (ii) implies (iii).
Conversely, if $I$ has a prime $p$ with multiplicity $m\geq 2$ then $\mu(I_p)=\mu(p_p^m)=m+1\geq 3$. Therefore, (iii) implies (ii).

Finally, quite generally for any codimension $2$ perfect ideal $I$ in a Cohen--Macaulay ring of dimension $3$, (iii) and (iv) are equivalent.
\qed

\medskip

A major question in this universe is how broad is the class of ideals of fat points allowing for an application of Corollary~\ref{red_at_most_2}.
Here we contend ourselves with a few special results.

The following notion was introduced in \cite[Section 4.2]{fat1} and more thoroughly discussed in \cite{fat2}.

\begin{Definition}\rm
	Let $\boldsymbol\mu=\{\mu_1,\ldots,\mu_r\}$ be a set of nonnegative integers  satisfying the following condition: there exists an integer $s\geq 2$ such that
	\begin{equation}\label{eqs-subhomaloidal}
	\sum_{i=1}^r \mu_i=3(s-1) \quad {\rm and}\quad \sum_{i=1}^r \mu_i^2=s(s-1).
	\end{equation}
	Note that the first of the above relations implies that $s$ is uniquely determined.
	We will say that $\boldsymbol\mu$ is a {\em  sub-homaloidal multiplicity set in degree $s$}.
\end{Definition}
We note that such an integer $s$ is necessarily odd.

\begin{Proposition}
Let $I\subset R=k[x,y,z]$ denote an ideal of plane fat points with a sub-homaloidal multiplicity set in degree $s\geq 3$ and such that $I=(I_s)$.Then:
\begin{enumerate}
	\item[{\rm (a)}] The image of the rational map defined by the linear system $I_s$ is a variety of degree $s/\mathfrak{d}$, where $\mathfrak{d}:=(k(I_s):k(R_s))$ is the field extension degree.
\end{enumerate}	
	Suppose in addition that the linear system $I_s$ has the expected dimension $(s+5)/2$ and that the special fiber $\mathcal{F}(I)$ is Cohen--Macaulay of dimension $3$. Then:
\begin{enumerate}	
	\item[{\rm (b)}] $\mathfrak{d}=1$, i.e., the rational map defined by $I_s$ is birational onto the image.
	\item[{\rm (c)}]  $r(I)=2$ and the coefficient of the highest term of the numerator of the Hilbert series of $\mathcal{F}(I)$ is $(s-1)/2$.
\end{enumerate}
\end{Proposition}
\demo (a) This is proved in \cite[Proposition 3.2]{fat2}.

(b)  Let $h:=1+h_1 t+h_2 t^2+\cdots$ stand for the numerator of the Hilbert series of $\mathcal{F}(I)$.
From (a) and the assumption on the expected dimension of $I_s$ one has: $$s/\mathfrak{d}=1+((s+5)/2-3)+\sum_{t\geq 2}h_t= (s+1)/2+\sum_{t\geq 2}h_t.$$ 
Since $\mathcal{F}(I)$ is Cohen--Macaulay, $\sum_{t\geq 2}h_t\geq 0$.
This is impossible unless $\mathfrak{d}=1$.

(c)  Since the rational map defined by $I_s$ is birational onto the image, the latter is a rational variety. By hypothesis, it is arithmetically Cohen--Macaulay, hence condition (iii) of Corollary~\ref{red_at_most_2} is satisfied (see, e.g., \cite[The proof of Corollary 2.6]{GGP}), thus implying that $r(I)\leq 2$.

In the notation of the proof of (b) we then obtain $h_t=0$ for $t\geq 3$.
This gives $h_2=s-(s+1)/2=(s-1)/2\geq 1$ and, in particular, $r(I)=2$.
\qed

\begin{Remark}\rm
Regarding the proof of item (c) in the above proposition, the reason to draw on Corollary~\ref{red_at_most_2}, instead the more direct Theorem~\ref{red_no2}, is that $I$ is not linearly presented. In \cite{fat2} an appropriate quadratic transform of $I$ is introduced that keeps the essential properties of the latter and, moreover, is a linearly presented ideal - for this Theorem~\ref{red_no2} could safely apply. 
\end{Remark}

\begin{Question}\rm
We note that the integer $(s-1)/2$ is a natural barrier for fat multiplicities of subhomaloidal types as explained in \cite{fat2}. How often does  $h_2=(s-1)/2$ hold while relaxing the second set of hypotheses in the above proposition?
\end{Question}

\subsection{Reciprocal ideals of hyperplane arrangements}

For the basic material of this part we refer to \cite{blup} and the references thereof.

Let $\mathcal{A}=\{H_1,\ldots,H_n\}\subset\mathbb P^{d-1}$ be a central hyperplane arrangement of size $n$ and rank $d$. Here $H_i=\ker(\ell_i),i=1,\ldots,n$, where each $\ell_i$ is a linear form in $R:=k[x_1,\ldots,x_d]$ and $k[\ell_1,\ldots,\ell_n]=R$.
We also assume that any two of these forms are linearly independent. From the algebraic viewpoint, there is a natural emphasis on the linear forms $\ell_i$ and the associated ideal theoretic notions.

With the above preamble, we consider the $(n-1)$-fold products of the linear forms $\ell_1,\ldots,\ell_n$:
\begin{equation}\label{products}
I:= (\ell_1\cdots \hat{\ell_i}\cdots \ell_n\,|\, 1\leq i\leq n).
\end{equation}

As it turns out, $I$ is a linearly presented codimension $2$ ideal with syzygy matrix having the shape
	\begin{equation}\label{phi4arrgs}
\varphi=\left(
\begin{array}{rrrr}
\ell_1&&&\\
-\ell_2&\ell_2&&\\
&-\ell_3&\ddots&\\
&&\ddots&\ell_{n-1}\\
&&&-\ell_n
\end{array}
\right)
	\end{equation}
where the blank slots are null entries.

The special fiber $\mathcal{F}(I)$ of this ideal is the so-called {\em Orlik--Terao algebra}. It is known to have dimension $d=\dim R$ (\cite{Te}) and to be Cohen--Macaulay (\cite{PrSp}). 

The following result has been proved in \cite[Lemma 3.1]{Schenck}. For the reader's convenience we give the proof  adapting to the present language.

\begin{Lemma}
	\label{Arr_is_fat}
If $I$ is as in {\rm (\ref{products})} then it is an ideal of fat points.	
\end{Lemma}
\demo 
Letting $\{\ell_1,\ldots,\ell_n\}$ denote the arrangement hyperplanes equations, where no two are proportional, we easily see that any associated prime $p$ of $R/I$ is generated by two such linear forms.
Say, $p=(\ell_1,\ell_2)$ without loss of generality.
Now consider all circuits of order $3$ involving $\ell_1,\ell_2$.
Let $c_p$ denote the number of these circuits and consider the $c_p$ additional linear forms intervening in these circuits, one for each such circuit. Without loss of generality, assume the additional forms are $\ell_3,\ldots,\ell_{t}$, where we have set $t:=c_p+2$.
Focusing on the sub-arrangement $\{\ell_1,\ell_2,\ldots,\ell_{t}\}$, we consider the ideal $J(p)$ of $R$ generated by its $(t-1)$-fold products.
Clearly, $J\subset p^{t-1}$. But since these ideals are generated in the same degree and the products generating $J(p)$ are $k$-linearly independent, we conclude that  $J(p)= p^{t-1}$.

Let $J:=\bigcap_{p\in {\rm Min} R/I} J(p)$.
Clearly, $I_p=J(p)_p$ for every $p\in {\rm Min} R/I$.
It follows that $I=J$, as contended.
\qed

\medskip

Assuming this result we prove that, if $\dim R=3$ then the  Rees algebra $\mathcal{R}(I)$ of $I$ is Cohen--Macaulay. This is a particular case of the general fact that $\mathcal{R}(I)$ is Cohen--Macaulay in any dimension, as proved in 
\cite{blup} by different methods.

\begin{Corollary}\label{Rees_of_Reciprocal_is_CM_dim3}
	Let $\dim R=3$ and let $I$ be as in {\rm (\ref{products})}.
	Then the Rees algebra $\mathcal{R}(I)$ is Cohen--Macaulay.
\end{Corollary}
\demo As remarked above, $\mathcal{F}(I)$ is Cohen--Macaulay.
By Theorem~\ref{red_no2}, it suffices to show that $r(I_p)\leq 1$ for every associated prime $p$ of $R/I$. Now $I$ is an ideal of fat points (in fact, in arbitrary dimension) by Lemma~\ref{Arr_is_fat}. For such an ideal, in dimension $3$, $I_p$ is generated by the power of a complete intersection in a $2$-dimensional ring. The latter is well-known to have reduction number $1$, with reduction generated by the two pure powers. 
\qed

A degenerate case is as follows:

\begin{Proposition}
Let $\dim R=3$ and let $I$ be as in {\rm (\ref{products})}.
The following are equivalent:
	\begin{itemize}
		\item[{\rm (i)}] Up to a change of variables in $R=k[x,y,z]$, the arrangement has the form
		$$\ell_1=x,\ell_2=y+c_2z,\ell_3=y+c_3z,\ldots,\ell_{n}=y+c_{n}z,$$
		where the coefficients $c_j\in k$ are mutually distinct.	
		\item[{\rm (ii)}] $u=1$one has
		\item[{\rm (iii)}] As a fat ideal, $I$ has the shape 
		$$I=q^{n-2}\cap\left( \bigcap _{p\in {\rm Min} R/I\setminus \{q\}} p\right),$$
		for a uniquely defined minimal prime $q$ of $R/I$.
	    \item[{\rm (iv)}] $r(I)\leq 1$	
	\end{itemize}   
\end{Proposition}
\demo 
(i) $\Rightarrow$ (ii)
The shape of $\phi$ as in (\ref{phi4arrgs}) becomes
$$\varphi=\left(\begin{matrix}
x&&&\\
-(y+c_2z)&y+c_2z&&\\
&-(y+c_3z)&\ddots&\\
&&\ddots&y+c_{n-1}z\\
&&&-(y+c_{n}z)
\end{matrix}
\right)$$
where the coefficients $c_j\in k$ are mutually distinct and the blank slots are null entries.
It is obvious that $I_1(\phi)=3$ and  $I_2(\phi)\subset (y,z)$. Therefore, $u=u(\phi)=1$ by definition.

(ii) $\Leftrightarrow$ (iii)
This follows from Proposition~\ref{chaos_4_fat}.

(ii) $\Rightarrow$ (iv)
This reads off Corollary~\ref{Lan1} (c).

(iv) $\Rightarrow$ (i) First, some preliminaries. We know that the special fiber $\mathcal{F}(I)$ of $I$ is Cohen--Macaulay, hence the assumption $r(I)\leq 1$ implies by the proof of Proposition~\ref{CMfiber} that $\mathcal{F}(I)$ is a rational normal scroll or a cone over a Veronese.
The homogeneous defining ideal of the latter is generated by the $2$-minors of a $2\times m$ matrix with the expected codimension $m-2+1=m-1$.
But since the dimension of $\mathcal{F}(I)$ is $\ell(I)=3$ then its codimension is $n-3$, and hence $m=n-2$.
Then
$$\mu(I^2)=\dim_k(I^2)_2=\dim_k \mathcal{F}(I)_2={n+1 \choose 2}-{n-2 \choose 2}=3(n-1).$$
From the other side, by Lemma~\ref{Arr_is_fat},  $I$ is an ideal of fat points. 
Letting $m_1,\ldots,m_r$ denote its defining virtual multiplicities, one has:
\begin{eqnarray*}\nonumber
\mu(I^2)&=&\dim_k I^2_{2(n-1)}\\  \nonumber
&=&\dim_k I^{(2)}_{2(n-1)} \quad\mbox{(by \cite[Corollary 1.4]{GGP})}\\ \nonumber
 &=&{2n\choose2}-\sum{2m_i+1\choose2} \quad\mbox{(by \cite[Corollary 1.4]{GGP})}\\ \nonumber
&=&n(2n-1)-2\sum_{i=1}^rm_{i}^2-\sum_{i=1}^rm_i.
\end{eqnarray*}
Therefore,
\begin{equation}\label{eq1}
2\sum_{i=1}^rm_i^2+\sum_{i=1}^rm_i=n(2n-1)-3(n-1)=2n^2-4n+3.
\end{equation}
Further,
$$n=\mu(I)=\dim_kI_{n-1}={n+1\choose 2}-\sum_{i=1}^r{m_i+1\choose2}=\frac{(n+1)n}{2}-\frac{\sum_{i=1}^rm_{i}^2+\sum_{i=1}^rm_i}{2}$$
It follows that  
\begin{equation}\label{eq2}
\sum_{i=1}^rm_{i}^2+\sum_{i=1}^rm_i=n^2-n
\end{equation}
From \eqref{eq1} e \eqref{eq2} we get the following formulas:
\begin{equation}\label{somadasmultiplicidades}
\sum_{i=1}^rm_{i}^2=n^2-3n+3\quad\mbox{e}\quad \sum_{i=1}^rm_i=2n-3.
\end{equation}
Next, we look at the nature of the arrangement.
Let, say, $\mathcal A'=\{\ell_1,\ldots,\ell_{s}\}$ stand for the largest subset of $\mathcal A=\{\ell_1,\ldots,\ell_n\}$ spanning a $k$-linear subspace of dimension $2.$ Clearly, $2\leq s.$ Also, $s<n$ since $\hht I_1(\phi)=3.$ 
If $s=n-1$ then, up to a change of variables, $\mathcal A$ would have the format stated in (i).
Thus, if we deny (i), it must be the case that $s\leq n-2$.
Let us argue that this leads to a contradiction.

We divide the argument in two cases:

\smallskip

\noindent{\bf First case:} $s>2.$

\smallskip

In this case, set $p:=(\ell_1,\ell_2)$, one of the minimal primes of $I$. 
The virtual multiplicity of $I$ at $p$ is $s-1\geq 2$ (see the proof of Lemma~\ref{Arr_is_fat}).
To get a grip on the other virtual multiplicities, we list a set of minimal primes of $I$ based on the index set $1\leq i\leq s$. 

For this, fixing any such $i$, we build a stratified partition $P(i)$ of the set $\mathcal A\setminus \mathcal A'$ in the following way: take some  $\ell\in \mathcal A\setminus \mathcal A'$ and define $P(i)_1$ to be the set consisting of $\ell$ and all elements of $\mathcal A\setminus \mathcal A'$ forming a circuit of order $3$ with $\ell_i$ and $\ell;$ take some $\ell\in (\mathcal A\setminus \mathcal A')\setminus P(i)_1$ and define $P(i)_2$ to be the set consisting of $\ell$ and all elements of $\mathcal A\setminus \mathcal A'$ forming a circuit of order $3$ with $\ell_i$ and $\ell;$ and so on so forth. 

Let $t(i)$ denote the cardinality of the partition $P(i)$.
For each $1\leq j\leq t(i)$, let $p_{i,j}:=(\ell_i, \ell_j)$, where $\ell_j$ is the element chosen at each step of the stratification $P(i)$. 
Clearly, together with $p=(\ell_1,\ell_2)$ these prime ideals are mutually distinct minimal primes of $I$.
Moreover, by construction, the virtual multiplicity of $I$ at $p_{i,j}$ is given by the cardinality, say, $m_{i,j}$, of the stratum $P(i)_j$  (again as in the proof of Lemma~\ref{Arr_is_fat}). 

It is now ripe time for counting.
One has $m_{i1}+\cdots+m_{is_i}=n-s$, hence
\begin{equation*}
\sum_{i=1}^rm_i\geq s-1+s(n-s)\geq s(n-s)+2.
\end{equation*}
Therefore,
\begin{eqnarray*}
\sum_{i=1}^rm_i-(2n-3)&\geq& s(n-s)-2n+5\nonumber\\
&=&n(s-2)-s^2+5\nonumber\\
&\geq& (s+2)(s-2)-s^2+5 \quad (\mbox{since}\;s\leq n-2)\nonumber\\
&=& 1,
\end{eqnarray*}
thus showing that $\sum_{i=1}^rm_i>2n-3.$ 
However, this contradicts the second equality in \eqref{somadasmultiplicidades}. 

\medskip

\noindent{\bf Second case:} $s=2.$

\medskip

This case is rather easy. Indeed, all combinations $(\ell_i,\ell_j)$ of mutually distinct elements of $\mathcal A$ yield ${n\choose2}$ distinct minimal primes of $I$.
This entails $\sum_{i=1}^rm_i\geq {n\choose2}>2n-3,$ also contradicting  the second equality in \eqref{somadasmultiplicidades}. 


\qed

\subsection{Linearly presented monomial ideals} \label{monomial}

As usual, we assume that $k$ is an algebraically closed field.

A linearly presented codimension $2$ primary ideal  $I\subset R=k[x,y,z]$ minimally generated by $n$ monomials is conjugate to the power $(y,z)^{n-1}$.


The corresponding zyzygy matrix is
\begin{equation}
\left(\begin{array}{ccccccccccc}
z&0&0&\ldots&0&0&0\\
-y&z&0&\ldots&0&0&0\\
0&-y&z&\ldots&0&0&0\\
\vdots&\vdots&\vdots&\ddots&\vdots&\vdots&\vdots\\
0&0&0&\ldots&z&0&0\\
0&0&0&\ldots&-y&z&0\\
0&0&0&\ldots&0&-y&z\\
0&0&0&\ldots&0&0&-y
\end{array}\right)
\end{equation}

Dealing with non-primary such ideals is somewhat more delicate.
Here is the ``next'' case:

\begin{Proposition}\label{canonical_two-primes} Let $I\subset R=k[x,y,z]$ denote a linearly presented codimension $2$ perfect ideal minimally generated by $n$ monomials of degree $n-1$. If $R/I$ has exactly two minimal primes then, up to conjugation, the syzygy matrix of $I$
has the form
	\begin{equation}\label{canonical_monomial}
\varphi=\left(\begin{array}{ccccccccccc}
	z&0&0&\ldots&0&0&0\\
	-c_1&z&0&\ldots&0&0&0\\
	0&-c_2&z&\ldots&0&0&0\\
	\vdots&\vdots&\vdots&\ddots&\vdots&\vdots&\vdots\\
	0&0&0&\ldots&z&0&0\\
	0&0&0&\ldots&-c_{n-3}&z&0\\
	0&0&0&\ldots&0&-c_{n-2}&z\\
	0&0&0&\ldots&0&0&-c_{n-1}
	\end{array}\right)
		\end{equation}
	where $c_i\in\{x,y\}$ for each $1\leq i\leq n-1.$
\end{Proposition}
\demo Up to conjugation -- in fact, up to change of variables in $R$ -- we may assume that the minimal primes of $R/I$ are $(x,z)$ and $(y,z)$, hence $\sqrt{I}=(xy,z)$. Therefore, a power of $z$ lies in $I$ and since $I$ is generated by monomials of degree $n-1$, it follows that $z^{n-1}\in I$.
By the same token, $I$ has a monomial generator of the form $x^ay^b$, with $a+b=n-1$.

We claim that a remaining set of $n-2$ monomial generators of $I$ is 
\begin{equation}\label{other_mon_gens}
\{x^{a_i}y^{b_i}z^i, \, 1\leq i \leq n-2, \,a_i+b_i=n-1-i\}
\end{equation}
In other words, no two minimal monomial generators have in their support the same power of $z$. Indeed, assuming this would be possible there would be a syzygy $\mathfrak{s}$ of $I$ with only two nonzero coordinates which are monomials in $x,y$ only.
Since one is assuming that the syzygy matrix $\phi$ of $I$ is linear, $\mathfrak{s}$ would have to be a combination of linear syzygies. Now, because $I$ is minimally generated by monomials, a generating linear syzygy can be assumed to have only two nonzero coordinates and these are then necessarily different signed variables. On the other hand, the generator $z^{n-1}$ being a maximal minor of $\phi$ implies that every column of $\phi$ has an entry $\pm z$ in exactly one row except one. It follows that $\mathfrak{s}$ cannot be a combination of the syzygies represented by the columns of $\phi$.

Having granted (\ref{other_mon_gens}), let us write the complete set of monomial generators thus obtained in the following ordering
\begin{equation}\label{gens_final_form}
\{x^{a_0}y^{n-1-a_0},\, x^{a_1}y^{n-2-a_1}z,\,\ldots, x^{a_{n-2}}y^{1-a_{n-2}}z^{n-2},\, z^{n-1} \}, a_i\geq 0.
\end{equation}
One immediately sees that a generating linear syzygy between these generators can only involve contiguous ones in the above ordering because of the sequential powers of $z$.
Moreover, its nonzero coordinates are $z$ and necessarily a signed $x$ or $y$. 
Therefore, $I$ admit at least those many syzygies that appear in the matrix (\ref{canonical_monomial}). But this partial matrix has rank $n-1$, hence must be a full syzygy matrix of $I$.
\qed

\medskip

We informally call the sequence of entries $c_1\,c_2\, \dots\, c_{n-1}$ the {\em basic entry sequence} of the ideal $I$ in the above assumptions.

We next give the structural properties of two remarkable cases of the basic entry sequence.

\begin{Theorem}\label{alternating}
Let the assumptions and notation be those of {\rm Proposition~\ref{canonical_two-primes}} and let the basic entry sequence be an alternating sequence of $x$ and $y$, i.e., of the shape $x\,y\,x\,y\,x\,y\dots$. Then:
\begin{enumerate}
	\item[{\rm (a)}] The special fiber of the ideal $I$ is Cohen--Macaulay and normal
	\item[{\rm (b)}] The Rees algebra of $I$ is Cohen--Macaulay and normal.
\end{enumerate}
\end{Theorem}
\demo We assume that the sequence has the shape $x\,y\,x\,y\,\dots\,x\,y\,x\,y$, i.e., that $n-1$ is even - the case where $n-1$ is odd is handled similarly with minor changes.
Denoting $\{f_0,f_1,\ldots,f_{n-1}\}$ the generators in (\ref{gens_final_form}), one has
$$f_i=\left\{\begin{array}{clccc}
x^{(n-1)/2}y^{(n-1)/2},& \mbox{if}\;i=0\\
(f_{i-1}/x)z,& \mbox{if}\; i\; \mbox{is odd} \\
(f_{i-1}/y)z,&\mbox{if}\;i\;\mbox{is even}
\end{array}\right.$$
(a) With this format it is straightforward to verify that the homogeneous defining ideal of the special fiber of $I$, which is the kernel of the $k$-algebra map
$$k[t_0,t_1,\ldots,t_{n-1}]\surjects k[f_0,f_1,\ldots,f_{n-1}], \, t_i\mapsto f_i,$$
contains the $2$-minors of the $2$-step catalecticant
$$C=\left(\begin{array}{ccccc}
t_0&t_1&t_2&\dots&t_{n-3}\\
t_2&t_3&t_4&\dots&t_{n-1}\end{array}\right).$$
It is well-known that $k[t_0,t_1,\ldots,t_{n-1}]/I_2(C)$ is a Cohen--Macaulay domain of codimension $n-3$. Therefore, this algebra is the defining ideal of $\mathcal{F}(I)$.

Clearly, this is normal since it is the homogeneous coordinate ring of a rational normal scroll (being catalecticant).

\smallskip

(b) We have seen in the proof of  (a)  that the homogeneous defining ideal of the special fiber of $I$ is the $2$-minors of the matrix
$$C=\left(\begin{array}{ccccccc}
t_0&t_1&t_2&\ldots&t_{n-3}\\
t_2&t_3&t_4&\ldots&t_{n-1}
\end{array}\right);$$
in particular, the fiber is Cohen--Macaulay and $r(I)\leq 1$.

We first show that the Rees algebra $\mathcal{R}_R(I)$ is Cohen--Macaulay. Since the fiber $\mathcal{F}(I)$ is Cohen--Macaulay, it suffices by Theorem~\ref{red_no2} to prove that $r(I_p)\leq 1$ and $r(I_q)\leq 1$ where $p=(x,z)$ and $q=(y,z).$ 
For this, we note that  
$$I=(x^my^m,x^{m-1}y^mz,\ldots,x^{i}y^iz^{n-1-2i},x^{i-1}y^iz^{n-2i+1},\ldots,xyz^{n-3},yz^{n-2},z^{n-1}),$$
with $m:=(n-1)/2.$ 
Then
$$
I_p=(x^m,x^{m-1}z,x^{m-2}z^3,\ldots,xz^{n-4},z^{n-2})\quad\mbox{and}\quad I_q=(y^m,y^{m-1}z^2,\ldots,yz^{n-3},z^{n-1})
$$
It follows that from \cite[Proposition 2.15]{GSWR} that $I_p$ and $I_q$ are normal ideals of $R_p$ and $R_q$, respectively. Thus, the Rees algebras $\mathcal{R}_{R_p}(I_p)$ and $\mathcal{R}_{R_q}(I_q)$ are Cohen--Macaulay. Therefore, $r(I_p)\leq 1$ and $r(I_q)\leq 1.$ 

We next prove that $\mathcal{R}_R(I)$ is normal, by showing that it satisfies Serre's condition $(R_1)$.
Let $\mathcal J$ stand for its presentation ideal on $R[\tt]$.
It contains the ideal $\mathcal I:=(\tt\cdot\phi,I_2(C))$ and the generators of the latter can be taken to be a subset of a full set of generators of $\mathcal J$.
Accordingly, the Jacobian matrix $\Theta(\mathcal I)$ of the generators of $\mathcal I$ is a submatrix of the Jacobian matrix $\Theta(\mathcal J)$ of the generators of $\mathcal J$, with same number of columns.
This gives $I_{n-1}(\Theta(\mathcal I))\subset I_{n-1}(\Theta(\mathcal J))$ and hence
$$(I_{n-1}(\Theta(\mathcal I), \mathcal I)\subset (I_{n-1}(\Theta(\mathcal J), \mathcal J),
$$
so it suffices to bound below the codimension of the leftmost ideal.
But the latter contains the ideal $(I_2(B)\cdot I_{n-3}(\Theta'), I_2(C), I)$
since 
$$\Theta(\mathcal I)
=\left(\begin{array}{ccccc}
B&0\\
\phi&\Theta'
\end{array}\right)
$$
where $\Theta'$ denotes the Jacobian matrix of the $2$-minors of $C$ and $B$ denotes the Jacobian dual of $\phi$.

Thus far the discussion has been kept on a fairly general level.
Now it is time to enter the particulars of the alternating case via the special shape of $B$.
Namely, the transpose of $B$ has the shape
$$
\left(\begin{array}{ccc}
-t_1&0&t_0\\
0&-t_2&t_1\\
-t_3&0&t_2\\
0&-t_4&t_3\\
\vdots&\vdots&\vdots\\
0&-t_{n-1}& t_{n-2}
\end{array}\right)
$$
Therefore, $(t_0t_{n-1},t_1^2, t_2^2,\ldots, t_{n-2}^2)\subset I_2(B)$, implying that the codimension of $(I_2(B),I)$ is at least $(n-1)+2=n+1$.
As to the other ``component'' $(I_{n-3}(\Theta'), I_2(C), I)$, one is through since $(I_{n-3}(\Theta'), I_2(C))$ is the Jacobian ideal of the normal domain $k[\tt]/I_2(C)$ of codimension $n-1$.
\qed

\begin{Conjecture} In the case above of an alternating basic entry sequence the Rees algebra of the ideal $I$ is  of fiber type.
\end{Conjecture}
The second remarkable case is as follows.

\begin{Theorem}\label{separating}
Let the assumptions and notation be those of {\rm Proposition~\ref{canonical_two-primes}} and let the basic entry sequence be a separating sequence of $x$ and $y$, i.e., of the shape $x x\,\dots\,x\,y\,y\,\dots\,y$. Then:
\begin{enumerate}
	\item[{\rm (a)}] The special fiber of the ideal $I$ is Cohen--Macaulay
	\item[{\rm (b)}] The Rees algebra of $I$ is Cohen--Macaulay and of fiber type.
\end{enumerate}
\end{Theorem}
\demo
(a) The argument in this case is more involved and it will appeal to a ``deformation'' argument.

Let $r$ (respectively, $s$) denote the number of $x$'s (respectively, the number of $y$'s) in the sequence $x x\,\dots\,x\,y\,y\,\dots\,y$.
By a similar token as in the previous case, $I$ is generated by the following $n-1$ forms
$$\left\{\begin{array}{clccc}
x^{r}y^{s},& \mbox{if}\;i=0\\
x^{r-i}y^sz^i,& \mbox{if}\; 1\leq i\leq r \\
y^{n-1-i}z^i,&\mbox{if}\;r+1\leq i\leq n-1
\end{array}\right.$$
Consider the following Hankel matrices
$$H_1=\left(\begin{array}{cccccccc}
t_0&t_1&\ldots&t_{r-1}\\
t_1&t_2&\ldots&t_{r}\end{array}\right)
\quad\mbox{and}\quad
H_2=\left(\begin{array}{cccccccc}
t_{r}&t_{r+1}&\ldots&t_{n-2}\\
t_{r+1}&t_{r+2}&\ldots&t_{n-1}\end{array}\right).
$$

{\sc Claim:} The ring $k[t_0,t_1,\ldots,t_{n-1}]/P$ is a Cohen--Macaulay domain of dimension $3$, where $P:=(I_2(H_1),I_{2}(H_2))k[t_0,t_1,\ldots,t_{n-1}]$.

\smallskip

We first show this ring is Cohen--Macaulay.
For this deform to $k[t_0,\ldots,t_{n-1},T],$ with $T$ a new variable, and consider the following Hankel matrix

$$H_2':=\left(\begin{array}{cccccccc}
T&t_{r+1}&\ldots&t_{n-2}\\
t_{r+1}&t_{r+2}&\ldots&t_{n-1}\end{array}\right).
$$
Setting $P':=(I_2(H_1),I_{2}(H_2'))\subset k[t_0,\ldots,t_{n-1},T]$, since $k$ is algebraically closed, we have that the $k$-algebra
$$S:=k[t_0,\ldots,t_{n-1},T]/P'\simeq k[t_0,\ldots,t_r]/I_2(H_1)\otimes_k k[t_{r+1},t_{r+2},\ldots t_{n-1},T]/I_2(H_2')
$$
is a Cohen--Macaulay domain of dimension $2+2=4$.
In particular, $T-t_r$ is regular on $S$, hence 
$$k[t_0,t_1,\ldots,t_{n-1}]/P\simeq S/(T-t_r)$$
is a Cohen--Macaulay ring of dimension $3$.

It remains to show that $P$ is a prime ideal.
For this we will show that $t_r$ is regular on $k[t_0,t_1,\ldots,t_{n-1}]/P$ and that the localization of the latter at $t_r$ is a domain.

Let $Q\subset k[t_0,t_1,\ldots,t_{n-1}]$ denote an associated prime of $k[t_0,t_1,\ldots,t_{n-1}]/P$ containing $t_r$.
We claim that necessarily $\{t_1,\ldots,t_{n-2}\}\subset Q,$ and hence $Q$ has codimension at least $n-2$; but this is impossible since $k[t_0,t_1,\ldots,t_{n-1}]/P$ is Cohen--Macaulay of codimension $n-3$.

Now, if $t_r\in Q$, since both minors 
$\det\left(\begin{array}{ccc}t_{r-2}&t_{r-1}\\t_{r-1}&t_r\end{array}\right)$ and $\det\left(\begin{array}{ccc}t_{r}&t_{r+1}\\t_{r+1}&t_{r+2}\end{array}\right)$
belong to $Q$, it follows that $t_{r-1}, t_{r+2}\in Q$.
But again, since $\det\left(\begin{array}{ccc}t_{r-3}&t_{r-2}\\t_{r-2}&t_{r-1}\end{array}\right)$ and $\det\left(\begin{array}{ccc}t_{r+1}&t_{r+2}\\t_{r+2}&t_{r+3}\end{array}\right)$ both belong to $Q$ it follows that $t_{r-2}, t_{r+2}\in Q$.
Continuing, eventually $\{t_1,\ldots,t_{n-2}\}\subset Q,$ as claimed. 

Finally, by inverting $t_r$ yields an isomorphism
$$\frac{k[t_0,\ldots,t_{n-1}][t_r^{-1}]}{Pk[t_0,\ldots,t_{n-1}][t_r^{-1}]}\simeq k[t_r,t_{r-1},t_{r+1}][t_{r}^{-1}].$$

\medskip

To close the proof, at the other end a direct verification gives that the ideals of $2$-minors $I_2(H_1)$ and $I_2(H_2)$ are contained in the defining ideal of the fiber $\mathcal{F}(I)$.
Therefore we conclude that  $k[t_1,\ldots,t_n]/P$ is the homogeneous coordinate ring of $\mathcal{F}(I)$.

\smallskip

(b) Consider the following scroll matrices
$$S_1=\left(\begin{array}{ccccccccc}
x&t_0&t_1&\ldots&t_{r-1}\\
-z&t_1&t_2&\ldots&t_{r}\end{array}\right)
\quad\mbox{and}\quad
S_2=\left(\begin{array}{ccccccccc}
y&t_{r}&t_{r+1}&\ldots&t_{n-2}\\
-z&t_{r+1}&t_{r+2}&\ldots&t_{n-1}\end{array}\right).
$$
and set $P:=(I_2(S_1),I_2(S_2))\subset k[x,y,z,t_0,\ldots,t_{n-1}]$.

We will show that $P$ coincides with a presentation ideal $\mathcal J$ of $\mathcal R (I)$.
For this it suffices to show that $P\subset \mathcal J$ and that $P$ is a prime ideal of dimension $4$. The inclusion is clear since the $2$-minors on each side are either defining relations of the symmetric algebra of $I$ or polynomial relations of $I$ (the latter by part (a)).
Moreover, once verified that $P=\mathcal J$, then it follows immediately that $\mathcal R(I)$ is of fiber type.

We proceed to show that $P$ is a Cohen--Macaulay prime ideal.

Let $T$ and $U$ new variables, and consider the following deformation of $S_2$ to  ``separate'' variables
$$S_2':=\left(\begin{array}{cccccccc}
y&U&t_{r+1}&\ldots&t_{n-2}\\
T&t_{r+1}&t_{r+2}&\ldots&t_{n-1}\end{array}\right).$$
Note that all three $S_1,S_2,S_2'$ define rational normal scroll surfaces in their respective projective ambients.

Set $P':=(I_2(S_1),I_{2}(S_2'))\subset k[x,y,z,t_0,\ldots,t_{n-1}, T,U]$.

{\sc Claim 1:} The ring $S := k[x,y,z,t_0,\ldots,t_{n-1}, T,U]/P'$ is a Cohen--Macaulay domain of dimension $6$.

\smallskip

Note that 
 $$k[t_0,\ldots,t_{n-1},x,y,z,T,U]/P'\simeq \frac{k[t_0,\ldots,t_r,x,z]}{I_2(S_1)} \otimes_k \frac{k[t_{r+1},t_{r+2},\ldots t_{n-1},y,T,U]}{I_2(S_2')}$$
is Cohen--Macaulay of dimension $3+3=6$.
Since $k$ is assumed to be algebraically closed, it is a domain.

Specializing back only one of the two deformation variables
$$S_2'':=\left(\begin{array}{cccccccc}
y&t_r&t_{r+1}&\ldots&t_{n-2}\\
T&t_{r+1}&t_{r+2}&\ldots&t_{n-1}\end{array}\right)$$
still give a matrix defining a rational normal scroll in the respective projective ambient.

Set $P'':=(I_2(S_1),I_{2}(S_2''))$.

{\sc Claim 2:} The ring $S' := k[x,y,z,t_0,t_1,\ldots,t_{n-1},T]/P''$ is a Cohen--Macaulay domain of dimension $5$.

\smallskip

Since this is an obvious specialization with $U-t_r$ regular on $S$, the assertion on the Cohen--Macaulayness and the dimension are immediate.

To show that $P''$ is a prime ideal we will argue that $t_r$ is regular on $S'$ and that the localization of the latter at $t_r$ is a domain.

Let $Q\subset k[x,y,z,t_0,t_1,\ldots,t_{n-1},T]$ denote an associated prime of $S'$ containing $t_r$. We claim that necessarily $\{t_1,\ldots,t_{n-2},zt_0,yt_{n-1}\}\subset Q,$ and hence $Q$ has codimension at least $n$; but this is impossible since $S'$ is Cohen--Macaulay of codimension $n-1$.

To see the claim, if $t_r\in Q$, since both minors $\det\left(\begin{array}{ccc}t_{r-2}&t_{r-1}\\t_{r-1}&t_r\end{array}\right)$ and $\det\left(\begin{array}{ccc}t_{r}&t_{r+1}\\t_{r+1}&t_{r+2}\end{array}\right)$
belong to $Q$, it follows that $t_{r-1}, t_{r+2}\in Q$. But again, since $\det\left(\begin{array}{ccc}t_{r-3}&t_{r-2}\\t_{r-2}&t_{r-1}\end{array}\right)$ and $\det\left(\begin{array}{ccc}t_{r+1}&t_{r+2}\\t_{r+2}&t_{r+3}\end{array}\right)$ both belong to $Q$ it follows that $t_{r-2}, t_{r+2}\in Q$. Continuing, eventually $\{t_1,\ldots,t_{n-2}\}\subset Q$. Moreover, both minors $\det\left(\begin{array}{ccc}x&t_0\\-z&t_1\end{array}\right)$ and $\det\left(\begin{array}{ccc}y&t_{n-2}\\T&t_{n-1}\end{array}\right)$ belong to $Q$, it follows that $zt_0,yt_{n-1}\in Q$. Therefore, $\{t_1,\ldots,t_{n-2},zt_0,yt_{n-1}\}\subset Q,$ as claimed.

Finally, inverting $t_r$ yields an isomorphism
$$\frac{k[x,y,z,t_0,t_1,\ldots,t_{n-1},T][t_r^{-1}]}{P'' k[x,y,z,t_0,t_1,\ldots,t_{n-1},T][t_r^{-1}]}\simeq k[y,z,t_r,t_{r-1},t_{r+1}][t_{r}^{-1}].$$

\medskip

{\sc Claim 3:} The ring $k[x,y,z,t_0,t_1,\ldots,t_{n-1}]/P$ is a Cohen--Macaulay domain of dimension $4$.

\smallskip

Again, this is a specialization with $T+z$ regular on $S'$, hence the assertions about the Cohen--Macaulayness and the dimension are immediate.

To see that $P$ is a prime ideal we follow the same path as above, namely, we will show that $t_r$ is regular on $S''$ and that the localization of the latter at $t_r$ is a domain.

Let $Q\subset k[t_0,t_1,\ldots,t_{n-1},x,y,z]$ denote an associated prime of $S''$ containing $t_r$. In the same way as above, we have that $\{t_1,\ldots,t_{n-2},zt_0,yt_{n-1}\}\subset Q,$ and hence $Q$ has codimension at least $n$; but this is impossible since $S''$ is Cohen--Macaulay of codimension $n-1$.

Finally, by inverting $t_r$ yields an isomorphism
$$\frac{k[x,y,z,t_0,\ldots,t_{n-1}][t_r^{-1}]}{P''k[x,y,z,t_0,\ldots,t_{n-1}][t_r^{-1}]}\simeq k[y,t_r,t_{r-1},t_{r+1}][t_{r}^{-1}].$$
\qed

\begin{Remark}\rm The special fiber of the separating case is not normal even when it is balanced (i.e., when the number of $x$ entries and of $y$ entries are equal).
	Basic entry sequences picked in random order of $x$ and $y$ has even less chance of having a normal special fiber since the latter may not even be Cohen--Macaulay.
	For example, the basic entry sequence $x\,y\,x\,x\,y\,y$, obtained by a single transposition of the alternating case, gives an ideal for which the special fiber is not Cohen--Macaulay. 
	In this example, $n=7$, while $u=3$, hence $n<2(u+1)$. A computation gives that $I_{u+1}(\phi)=I_4(\phi)$ has two minimal primes. This should be confronted with Proposition~\ref{chaos_and_minimal_primes} (b).
	
	The question remains as to whether there are other events in the present context where normality takes place besides the alternating (scroll) case.
\end{Remark}


\noindent {\bf Authors' addresses:}

\medskip

\noindent {\sc Andr\'e Doria},  Departamento de Matem\'atica, CCET\\ Universidade Federal de Sergipe\\
49100-000 S\~ao Cristov\~ao, Sergipe, Brazil\\
{\em e-mail}: avsdoria@ufs.br\\

\noindent {\sc Zaqueu Ramos},  Departamento de Matem\'atica, CCET\\ Universidade Federal de Sergipe\\
49100-000 S\~ao Cristov\~ao, Sergipe, Brazil\\
{\em e-mail}: zaqueu.ramos@gmail.com\\

\noindent {\sc Aron Simis},  Departamento de Matem\'atica, CCEN\\ Universidade Federal de Pernambuco\\ 50740-560 Recife, PE, Brazil\\
{\em e-mail}:  aron@dmat.ufpe.br

\end{document}